\documentclass[3p]{elsarticle}

\makeatletter
\def\ps@pprintTitle{%
    \let\@oddhead\@empty
    \let\@evenhead\@empty
    \let\@evenfoot\@oddfoot
    }
\makeatother

\def\text#1{\mbox{#1}}

\def\be{\begin{equation}}
\def\ee{\end{equation}}
\def\ba{\begin{array}}
\def\ea{\end{array}}
\def\bea{\begin{eqnarray}}
\def\eea{\end{eqnarray}}

\def\bt{\begin{tabular}}
\def\et{\end{tabular}}
\def\bc{\begin{center}}
\def\ec{\end{center}}

\def\illustration #1 by #2 (#3){
  \vbox to #2{
    \hrule width #1 height 0pt depth 0pt
    \vfill
    \special{illustration #3}}}

\def\scaledillustration #1 by #2 (#3 scaled #4){{
  \dimen0=#1 \dimen1=#2
  \divide\dimen0 by 1000 \multiply\dimen0 by #4
  \divide\dimen1 by 1000 \multiply\dimen1 by #4
  \illustration \dimen0 by \dimen1 (#3 scaled #4)}}

\def\scaledpicture #1 by #2 (#3 scaled #4){{
  \dimen0=#1 \dimen1=#2
  \divide\dimen0 by 1000 \multiply\dimen0 by #4
  \divide\dimen1 by 1000 \multiply\dimen1 by #4
  \picture \dimen0 by \dimen1 (#3 scaled #4)}
  }

\def\WidthBoarder{0pt}

\def\picture #1 by #2 (#3){
\vrule width \WidthBoarder height #2 depth 0pt
  \vbox to #2{
    \hrule width #1 height \WidthBoarder depth 0pt
    \vfill
    \special{picture #3} 
\hrule width #1 height \WidthBoarder depth 0pt
             }
\vrule width \WidthBoarder height #2 depth 0pt
                          }

\def\draft{0}

\def\EPSF#1#2#3#4#5{
\begin{figure}[p]
\vskip 0pt plus 0.5in
$$
\scaledillustration #3 by #4 (#2.epsf scaled #5)
$$
\if\draft1\caption{#1 (File: #2)}
\else
\caption{#1}
\fi
\vspace{0pt plus 0.5in}
\label{#2}
\end{figure}
}

\def\EPS#1#2#3{
\begin{figure}
\vskip 0pt plus 0.5in
$$
\BoxedEPSF{#2.eps scaled #3}
$$
\if\draft1\caption{#1 (File: #2)}
\else
\caption{#1}
\fi
\vspace{0pt plus 0.5in}
\label{#2}
\end{figure}
}

\def \text#1{\hbox{\rm #1}}

\def \ctrline#1{\hbox to \hsize{\hfil#1\hfil}}	

\def \newpage{\par \vfill \eject}

\def\smallarray#1\\#2\endsmallarray
    {{\def\\{\cr}
      \def\c{$\hfill\scriptstyle{####}\hfill$}
      \def\r{$\hfill\scriptstyle{####}$}
      \def\l{$\scriptstyle{####}\hfill$}
      \xdef\rowspecs{#1}
   \vcenter{\baselineskip 8pt\lineskip1pt
       \halign{\tabskip 0pt\rowspecs\cr#2\cr}}}}

\usepackage{graphics}
\usepackage{graphicx}
\usepackage{amssymb}
\usepackage{amsthm}
\usepackage{subfigure}
\usepackage{float}
\usepackage{amsmath}
\usepackage{amssymb}
\usepackage{mathrsfs}
\usepackage{relsize}
\usepackage{booktabs}
\usepackage{lineno}
\usepackage{wrapfig}
\usepackage[export]{adjustbox}
\graphicspath{{../Figures/}}
\usepackage[dvipsnames]{xcolor}
\usepackage{caption}
\biboptions{sort&compress}

\usepackage{pifont}
\usepackage{hhline}
\newcommand{\cmark}{\text{\ding{51}}}
\newcommand{\xmark}{\text{\ding{55}}}

\usepackage{hyperref}

\theoremstyle{remark}
\newtheorem{remark}{Remark}[section]

\journal{Computer Methods in Applied Mechanics and Engineering}

\begin{document}
\begin{frontmatter}
\title{Inference of Heterogeneous Material Properties via\\
 Infinite-Dimensional Integrated DIC}

\author[label1]{Joseph Kirchhoff\corref{cor1}}\ead{kirchhoff@my.utexas.edu}
\author[label2]{Dingcheng Luo}
\author[label2]{Thomas O'Leary-Roseberry}
\author[label1,label2]{Omar Ghattas}
\address[label1]{Walker Department of Mechanical Engineering, The University of Texas at Austin, 204 E. Dean Keaton Street, C2200, Austin, TX 78712, USA}
\cortext[cor1]{Corresponding author}
\address[label2]{Oden Institute for Computational Engineering and Sciences, The University of Texas at Austin, 201 E. 24th Street, C0200, Austin, TX 78712, USA}

\begin{abstract}

We present a scalable and efficient framework for the inference of spatially-varying parameters of continuum materials from image observations of their deformations.
Our goal is the nondestructive identification of arbitrary damage, defects, anomalies and inclusions without knowledge of their morphology or strength. Since these effects cannot be directly observed, we pose their identification as an inverse problem.
Our approach builds on integrated digital image correlation (IDIC, Besnard Hild, Roux, 2006), 
which poses the image registration and material inference as a monolithic inverse problem, thereby enforcing physical consistency of the image registration using the governing PDE.
Existing work on IDIC has focused on low-dimensional parameterizations of materials.
In order to accommodate the inference of heterogeneous material property fields that are formally infinite dimensional,
we present $\infty$-IDIC, a general formulation of the PDE-constrained coupled image registration and inversion posed directly in the function space setting. 
This leads to several mathematical and algorithmic challenges arising from the ill-posedness and high dimensionality of the inverse problem. 
To address ill-posedness, we consider various regularization schemes, namely $H^1(\Omega)$ and total variation for the inference of smooth and sharp features, respectively. 
To address the computational costs associated with the discretized problem, we use an efficient inexact-Newton CG framework for solving the regularized inverse problem. 
In numerical experiments, we demonstrate the ability of $\infty$-IDIC to characterize complex, spatially varying Lamé parameter fields of linear elastic and hyperelastic materials. 
Our method exhibits (i) the ability to recover fine-scale and sharp material features, (ii) mesh-independent convergence performance and hyperparameter selection, (iii) robustness to observational noise.

\end{abstract}


\end{frontmatter} 

\tableofcontents


\section{Introduction}

We present a novel formulation of inverse problems for spatially-varying parameter fields of continuum materials using images of their deformations.
In particular, we simultaneously pose the inversion and the image registration problems as a continuum limit in a function space setting. 
This formulation allows for the derivation of scalable and efficient \emph{discretization dimension-independent} algorithms for the inversion of heterogeneous material properties such as Lam\'{e} parameters while maintaining sharpness of features in their reconstruction, via the use of appropriate regularization methods. 
Our approach aims to address the imperative to characterize materials with spatially-varying properties. In particular we seek the non-destructive identification of arbitrary damage, defects, degradation, anomalies and inclusion without knowledge of their morphology or strength. 
The need for such methods is laid out in the vision for Materials Testing 2.0 (MT2.0) \cite{Pierron2023MaterialReview, Pierron2021TowardsMeasurements}, which calls for a paradigm shift in material characterization within research laboratories, particularly for materials like tissues, composites, welds, foams, etc., where full-field measurements should be used to identify spatially varying properties \cite{Pierron2021TowardsMeasurements}.

Since heterogeneous material defects are not directly observable, their identification must be posed as an inverse problem for arbitrary infinite-dimensional fields. In order to infer such infinite-dimensional fields, a significant amount of \emph{informative} observational data is required. Digital image correlation (DIC) derives displacement data via statistical correlation in differences in speckle patterned images \cite{Blaber2015Ncorr:Software}. However since these statistical algorithms do not obey the governing mechanical equilibrium equations, they may lead to cascading errors in the inversion. Integrated DIC (IDIC)  introduced by Besnard, Hild and Roux, overcame this issue by constraining the displacements to themselves follow the governing equations of the deformation process \cite{Besnard2006Finite-elementBands}. IDIC in its current formulation however inverts for material properties as homogeneous constants \emph{and not} heterogeneous spatial fields \cite{Kosin2024}. 

In this paper, we extend the formulation of IDIC to the function space setting, where we formulate the image registration and heterogeneous material inference problem as one monolithically coupled inverse problem. In this setting, we model the image data, and the material properties as infinite-dimensional spatial fields alongside the state variable (e.g., strain). Since the material properties are now formally infinite-dimensional fields, this leads to significant ill-posedness. In order to overcome this ill-posedness we consider various regularization models, such as $L^2(\Omega)$, $H^1(\Omega)$ Tikhonov regularization for smooth parameter reconstruction, and (primal-dual) total variation (TV) for reconstruction of sharp features. These modeling choices allow us to derive scalable and efficient Newton methods for the solution of these corresponding inverse problems. We term the methods presented herein Infinite-Dimensional IDIC, abbreviated as $\infty$--IDIC. 

We illustrate our approach through nondestructive evaluation problems for solids where one has observations of undeformed and deformed states of the material, which are painted speckle patterns.

In the numerical results we demonstrate $\infty$--IDIC considering linear elastic and neo-Hookean hyperelastic models. We demonstrate that $\infty$--IDIC can recover spatially-varying material properties with remarkable accuracy and sharpness; the corresponding stress fields can be post-processed for. This shows promise for predicting failure mechanisms of stress concentrations that arise from spatially varying material properties. Moreover, that our algorithms exhibit dimension-dependent convergence properties, and robustness to observational and loading noise. Our findings show that the information gained by the inversion depends on the loading conditions, which can be exploited by incorporating multiple experiments in a single $\infty$--IDIC inverse problem.

\subsection{Related works}

IDIC has been demonstrated using a wide range of PDEs such as phase-field fracture modeling \cite{Kosin2024}, elasto-plasticity \cite{Mathieu2015EstimationIntegrated-DIC,Neggers2019SimultaneousIdentification, Li2022Local-Micro-Zone-WiseJoints}, anisotropic linear elasticity \cite{Lindstrom2023IntegratedPlates}, to name a few. The current research in IDIC focuses on low-dimensional inversion, where the problem's ill-posedness may not be as severe as in high-dimensional spaces. However, often Tikhonov regularization is opted for or no regularization at all \cite{Neggers2019SimultaneousIdentification}. Rokos et al. address low-dimensional parameter estimation limitations by using IDIC to infer a homogenized modulus and then employing the principle of virtual work for heterogeneous behavior inference \cite{Rokos2023}. Such studies underscore the necessity for a direct approach to handle IDIC in high-dimensional function spaces. 

Unlike IDIC, which directly solves the inverse problem for material parameters, an alternative approach involves using image registration DIC software to acquire a displacement field, subsequently treating it as observational data. Initially devised to validate Finite Element Model (FEM) predictions against DIC observations \cite{Kavanaght1971FINITESOLIDSt}, a technique known as Finite Element Model Updating (FEMU), has evolved into the most commonly used method for material parameter inference from DIC data \cite{Sun2005FiniteMethod, Lenny2009, Waisman2010DetectionAlgorithms, Elouneg2021, Ereiz2022ReviewApplications, Fayad2023OnUpdating}. It was shown that FEMU and IDIC behave similarly in a simple experiment for inferring elasto-plastic parameters \cite{Mathieu2015EstimationIntegrated-DIC}. However, IDIC outperforms FEMU in challenging inverse problems, such as those involving image noise, complex loading conditions or experimental errors \cite{Ruybalid2016ComparisonDIC}. FEMU is useful when the DIC algorithm accurately captures the displacement field, yet any error in DIC image registration can compromise the material inversion \cite{Lehoucq2021TheCorrelation}. 

Existing IDIC works infer scalar, homogeneous material parameters, largely neglecting spatially varying fields, specifically heterogeneous materials \cite{Pierron2021TowardsMeasurements, Pierron2023MaterialReview}. There is limited work in FEMU that invert for spatially varying material properties for one or two stiff inclusions in soft bodies for tumor identification (elastography) \cite{skovoroda1999reconstructive,zhu2003finite, liu2005tomography, goenezen2011solution, mei2016estimating, mei2018comparative}. The existing work is reliant on an initial DIC algorithm to provide accurate displacement fields and do not appear to be mesh-independent. Scaling to high dimensional fields is a challenge we address in this paper for inversion of complex material fields. Infinite dimensional inverse problems, where the input parameters are \emph{inferred} from noisy observations of the data, are inherently ill-posed \cite{,engl1996regularization,vogel2002computational}. We solve this using an inverse problem formulation that is familiar to a vast literature. For more information we refer the reader to other works \cite{engl1996regularization,vogel2002computational, Tarantola2005, Stuart2010,Ghattas2021} for a more in depth treatment of the mathematical aspects of inverse problems.

\subsection{Contributions}

We present $\infty$--IDIC, an image based parameter estimation inverse problem established in function spaces, thereby enabling inversion of spatially varying material properties. We note that the formulation of $\infty$--IDIC is general and applies across various materials like tissues, composites, welds, and foams, improving material characterization through full-field measurements. Inferring detailed parameter and stress fields offers new possibilities, a notable advancement for both FEM validation and non-destructive evaluation. 

\begin{figure}[H]
    \centering
    \includegraphics[width=0.70\textwidth]{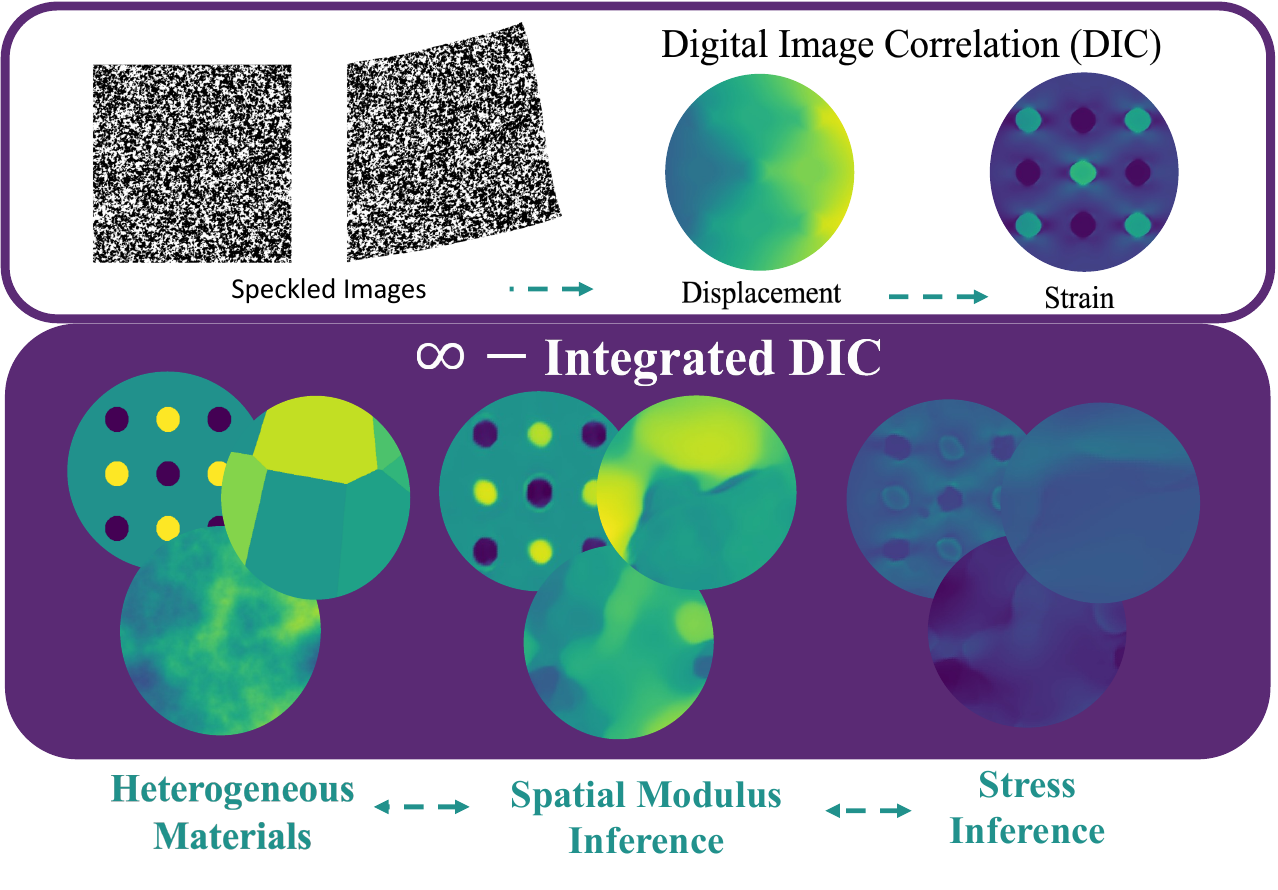}
    \caption{Our contributions include the development of an infinite-dimensional IDIC framework that enables the inversion of spatially varying parameter fields engineered for heterogeneous materials. Novel contributions are highlighted: handling heterogeneity, infinite dimensional framework, inference for high dimensional modulus fields, and stress fields.}
    \label{fig:IDIC_framework}
\end{figure}

Our overall approach of simultaneously posing the image registration and heterogeneous material inference problems can be extended to other fields such as in fluids (particle tracking velocimetry) \cite{Maas1993ParticleCoordinates, Malik1993ParticleTracking, Abdulwahab2020AApplications, Brevis2011IntegratingVelocimetry, Guezennec1994AlgorithmsVelocimetry, Schroder2022AnnualMechanics}, rheometry (micro-rheometry) \cite{Ahmadzadegan2023},  medical (image registration) \cite{Wyawahare2009ImageOverview, GHill2001PhysicsRegistration, Oliveira2014MedicalReview, Fu2020DeepReview, Sotiras2013DeformableSurvey, Rigaud2019DeformableEvaluation}, robotics (point-set registration) \cite{Sandhu2010PointDynamics, Ma2016Non-rigidStructures, Yuan2023Non-rigidChallenges, Maiseli2017RecentMethods, Min20233-DPerspective}, and even more broadly for optical flow \cite{Liu2015ComparisonImages, Wang2015AnFlow, Shah2021TraditionalInvestigation}. Figure \ref{fig:IDIC_framework} illustrates the key contributions of our work, by highlighting the novelty of handling heterogeneous material properties, spatial modulus inference and ultimately an accurate stress inference. The key contributions are summarized:

\begin{enumerate}
    \item Formulating $\infty$--IDIC in infinite-dimensional function spaces, enabling a dimension-independent algorithm to infer spatially varying parameter fields.
    \item Inverting for high-dimensional modulus fields captures heterogeneous material properties. 
    \item Post-processing for stress fields from the inferred modulus and displacement fields, enabling the identification of heterogeneity induced stress concentrations.
\end{enumerate}

We handle the inherent ill-posedness of the formulated inverse problem via mesh-independent regularization strategies, including $L^2(\Omega)$, $H^1(\Omega)$, primal-dual total variation that result in modeling benefits (e.g., sharpness, smoothing). This results in a regularized inverse problem which we choose to solve with an inexact Newton-CG method using a Gauss--Newton approximation of the Hessian. We demonstrate that the inferred modulus and displacement fields can be used to predict stress fields, enabling the identification of stress concentration-based failure mechanisms. In our numerical experiments, linear elastic and hyperelastic material models are considered; these results showcase the robustness of $\infty-$IDIC to noise, mesh independence, and boundary condition variations for complex material systems. We highlight our contributions relative to existing works in Table \ref{table:contributions}, and illustrate the capabilities of $\infty$--IDIC in Figure \ref{fig:IDIC_framework}. 
\begin{table}
\centering
{\renewcommand{\arraystretch}{1.3}
\begin{tabular}{|c|c|c|c|} 
\hline
 & \begin{tabular}{@{}c@{}}physics-based \\ registration\end{tabular} & heterogeneity & \begin{tabular}{@{}c@{}}mesh \\ independence\end{tabular}    \\ 
\hline
FEMU (scalar) \cite{Sun2005FiniteMethod, Lenny2009, Waisman2010DetectionAlgorithms, Mathieu2015EstimationIntegrated-DIC, Elouneg2021, Ereiz2022ReviewApplications, Fayad2023OnUpdating} &  {\color{BrickRed} \xmark} & {\color{BrickRed} \xmark} & {\color{BrickRed} \xmark}    \\ 
\hline
FEMU (field) \cite{skovoroda1999reconstructive,zhu2003finite, liu2005tomography, goenezen2011solution, mei2016estimating, mei2018comparative} &  {\color{BrickRed} \xmark} & {\color{ForestGreen} \cmark} & {\color{BrickRed} \xmark}    \\ 
\hline
IDIC (scalar) \cite{Besnard2006Finite-elementBands, Mathieu2015EstimationIntegrated-DIC,Neggers2019SimultaneousIdentification, Li2022Local-Micro-Zone-WiseJoints, Lindstrom2023IntegratedPlates, Kosin2024} & {\color{ForestGreen} \cmark}  & {\color{BrickRed} \xmark} & {\color{BrickRed} \xmark}    \\
\hline
$\infty$-IDIC (ours) & {\color{ForestGreen} \cmark}  & {\color{ForestGreen} \cmark}  & {\color{ForestGreen} \cmark}    \\
\hline
\end{tabular}
}
\caption{Summary of different DIC-based material property inverse methods.}
\label{table:contributions}
\end{table}

\subsection{Layout of the paper}
The rest of the paper is divided into two overarching sections: formulation of $\infty$--IDIC and numerical results. The formulation section begins with describing the infinite-dimensional inverse problem and its ill-posedness. This is treated using regularization methods, including $L^2(\Omega)$ Tikhonov regularization, $H^1(\Omega)$ Tikhonov smoothing regularization, and total variation to favor sharp reconstructions. We then present an inexact-Newton algorithm for solving the resulting PDE-constrained optimization problems, including the computation of the gradient and (approximation) Hessian actions via adjoint methods. The numerical results' section showcases the performance of the algorithm in inferring high-dimensional modulus fields. The results are presented for linear elasticity and hyperelasticity, demonstrating the robustness of the algorithm to noise, mesh independence, and boundary condition variations. The paper concludes with a discussion on the results and future work.


\section{The $\infty$--IDIC Formulation}

\begin{figure}[H]
    \centering
    \includegraphics[width=0.50\textwidth]{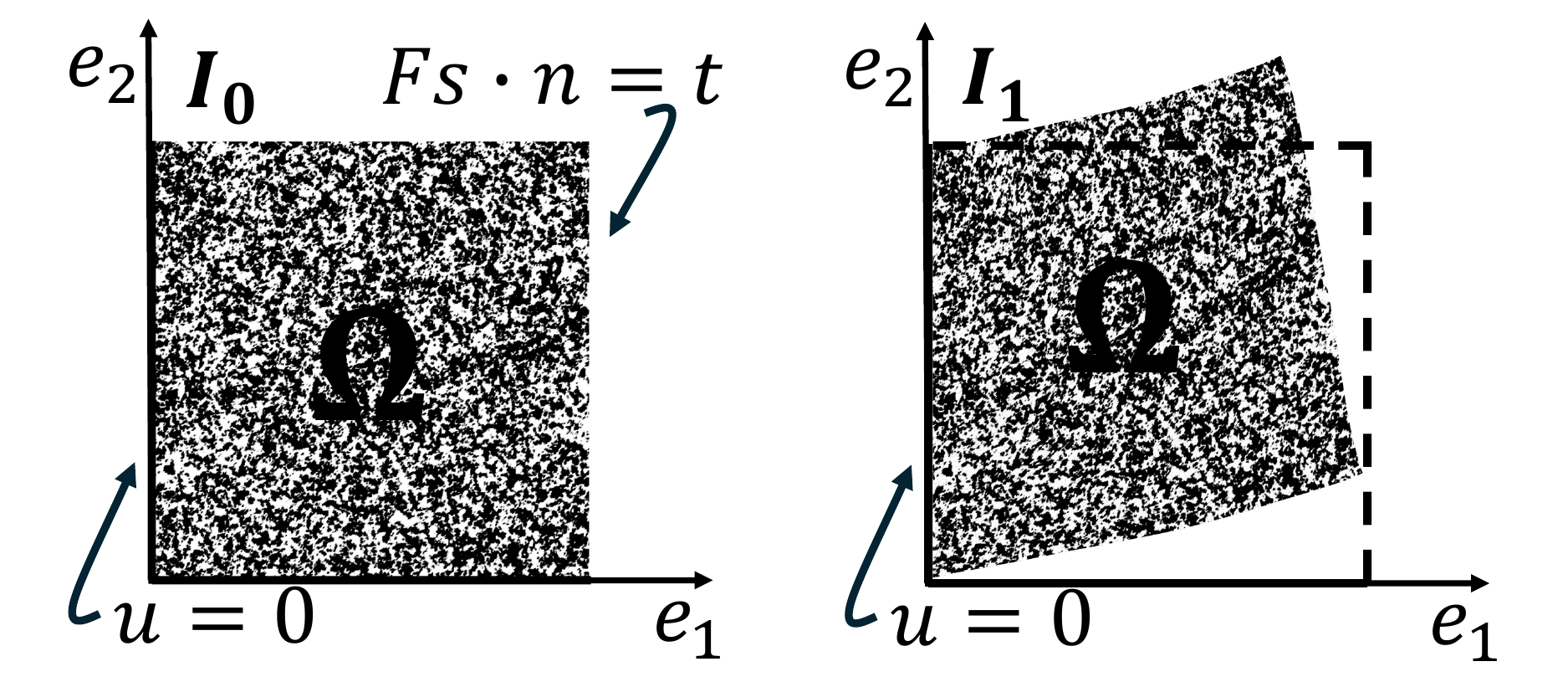}
    \caption{Schematic of the numerical setup for generating the reference and deformed states. The domain, $\Omega$, is speckled. A traction condition is applied to the right boundary, $\Gamma_R$ and a clamped condition is applied to the left boundary, $\Gamma_L$.}
    \label{fig:deformation_schematic}
\end{figure}

This section presents the formulation for $\infty$--IDIC: the simultaneously coupled infinite-dimensional integrated image registration and parameter identification problem from image data. We assume a given model for the underlying physics (e.g., linear elasticity and hyperelasticity), and assume that the observational data for the inverse problem consist of images of the undeformed and deformed configurations of the body of interest. We assume that the body is painted with a speckle pattern (speckled) to maximize the information about the point wise deformation state via contrast in pixel values (see Figure \ref{fig:deformation_schematic}). We then formulate a simultaneously coupled image registration and heterogeneous material property inverse problem, which constrains the admissible displacements to obey the governing conservation laws. This formulation thereby leads to physically consistent reconstruction of both the displacement and heterogeneous parameter fields. We follow the ``optimize-then-discretize'' (OTD) formalism, and derive all of our proposed methods for the continuum limit of the problems of interest in their corresponding function space setting \cite{becker2007optimal,braack2009optimal}. This allows us to mitigate the effects of discretization on the performance of the methods in practice.

We consider physical models governing the motion of deforming elastic bodies, formulated as partial differential equations (PDEs). 
Given a physical domain $\Omega \subset \mathbb{R}^d$ ($d = 2,3$) and boundary $\Gamma$ of the (undeformed) solid body, 
the PDE model defines the mapping from a spatially varying parameter field $m\in \mathcal{M}$ and the traction $t \in \mathcal{T}$ 
to the state variable $u \in \mathcal{U} := u_0 + \mathcal{V}$ representing the displacement field, where $\mathcal{M}$, $\mathcal{T}$, $\mathcal{U}$ are the function spaces for parameter, traction, and displacement fields, respectively.
A schematic of an example setup is shown in Figure \ref{fig:deformation_schematic}.
In particular, $\mathcal{M}$ and $\mathcal{T}$ 
are assumed to be Hilbert spaces such as the Sobolev spaces $H^k(\Omega)$ and $L^2(\Gamma)$.
The state space $\mathcal{U}$ is given as an affine shift of a Hilbert space $\mathcal{V}$ 
by a finite energy lift of any displacement boundary conditions, $u_0$, 
where $u_0$ satisfies the desired Dirichlet boundary conditions and $\|u_0\|_{\mathcal{V}} < \infty$. 
The PDE model can then be written abstractly as
\begin{equation} \label{eq:forward_pde}
	\text{PDE Model:} \qquad R(u,m,t) = 0, 
\end{equation}
where the PDE residual, $R: \mathcal{U} \times \mathcal{M} \times \mathcal{T} \rightarrow \mathcal{V}'$,
is a possibly nonlinear combination of differential operators,
and $\mathcal{V}'$ is the (topological) dual of $\mathcal{V}$. 
Alternatively, we can formulate the PDE problem in its weak form using $\mathcal{V}$ as the test space, 
\begin{equation} \label{eq:forward_pde_weak}
	\text{Find } u \in \mathcal{U} \text{ such that} \qquad r(u,m,t,v) = \langle R(u,m,t), v \rangle_{\mathcal{V}} = 0, \qquad \forall v \in \mathcal{V},
\end{equation}
where $r : \mathcal{U} \times \mathcal{M} \times \mathcal{T} \times \mathcal{V} \rightarrow \mathbb{R}$ 
is the weak residual, and $\langle \cdot, \cdot \rangle_{\mathcal{V}}$ denotes the duality pairing between $\mathcal{V}'$ and $\mathcal{V}$.
Note that $r(u,m,t,v)$ is linear with respect to the test function $v$. 
We will assume that the PDE problem is well-posed and admits solution operator $u = u(m,t)$ that is differentiable with respect to the PDE parameters.

In the simplest case, the observational data is given in the form of two images: $I_0$ and $I_1$, capturing the undeformed and deformed specimen, respectively, see Figure \ref{fig:deformation_schematic} for a schematic. This case generalizes to many images in the case of multiple deformation states to better inform the inversion. 
Mathematically, we consider images as functions of spatial coordinates that return the grayscale pixel value of the image at a particular point, 
where $0$ is black and $255$ is white. 
In particular images $I_0$ and $I_1$ are defined over an image domain $\Omega_I$, which is large enough to capture both the undeformed and deformed specimen. 
Moreover, we consider the space of images to be $\mathcal{I} = H^1(\Omega_I)$ such that they admit at least one spatial derivative.

The objective of the coupled inversion and registration is then to find the parameter field $m$ such that the predicted displacement yields
an image that matches the image of the deformed specimen.
This is formulated using the following minimization problem,
\begin{equation} \label{eq:idic_minimization_problem}
	\min_{m \in \mathcal{M}, u \in \mathcal{U}} \mathcal{J}(u,m; I_0, I_1) := \Phi(u; I_0,I_1) + \mathcal{R}(m)  \text{ subject to } R(u,m,t) = 0, \nonumber
\end{equation}
Here, the image misfit
\begin{equation} 
	\Phi(u; I_0,I_1) = \frac{1}{2}\int_{\Omega} |I_1(x+u(x)) - I_0(x)|^2 dx,
\end{equation}
compares the pixel values of the images, where $| \cdot |$ denotes the magnitude. 
Simply put, the misfit is the difference in the two images in function space where the displaced image is `pulled-back' by the displacement field. 

For the elliptic and parabolic PDEs that govern the deformation of solid bodies, 
the mapping from material parameter fields $m$ to the displacement $u$ tends to be smoothing, 
i.e., highly oscillatory perturbations of $m$ do not significantly affect $u$. 
Thus, the inversion of full heterogeneous parameter fields is inherently ill-posed.
To overcome this ill-posedness, we additionally introduce regularization on the parameter field, $\mathcal{R}(m)$, to the cost functional $\mathcal{J}$.

It is important to note that, unlike in the case of DIC, where the ill-posedness is in directly identifying the displacement field based on pixel correlations, the displacement field in our formulation is completely constrained by the underlying PDE solution, $u = u(m,t)$. 
Instead, the regularization is prescribed on the parameter field to overcome the ill-posedness arising from the desire to invert for heterogeneous material parameters.

\subsection{Handling the Ill-posedness of the Inverse Problem}  \label{sec:regularization}

As previously discussed, it is typical for the PDE solution operator to be insensitive to certain directions (subspaces) in the parameter space.
In such settings, observation data alone is not sufficient to identify the material parameters in these subspaces, 
and regularization is required to make the inverse problem well-posed. 
The choice of regularization $\mathcal{R}$ is based on the prior knowledge about the parameter field, 
and is tightly related to the choice of a prior distribution on the parameters in Bayesian inference.
Thus, the context should motivate the decision of which regularization to use. 

In this work, we consider three choices of regularization terms. 
In the simplest case, we can consider $L^2$ Tikhonov regularization to penalize the deviation of $m$ from a reference value $\bar{m} \in \mathcal{M}$,
\begin{equation}
	\mathcal{R}_{L^2}(m) = \frac{\gamma_{L^2}}{2} \int_{\Omega} |m(x) - \bar{m}(x)|^2 \, dx.
\end{equation}
where $\gamma_{L^2} > 0$ is a weighting factor.
For example, one can take $\bar{m}(x) \equiv m_{\text{nominal}}$, where $m_{\text{nominal}}$ is a nominal value for the material under consideration, maybe obtained from a material database.
In addition to directly regularizing the values of $m$, 
another option is to penalize steeper gradients of $m$ 
to impose a preference for smooth parameter fields. 
To this end, $H^1$ Tikhonov regularization uses an $L^2(\Omega)$ penalization on the gradient,
\begin{equation}
	\mathcal{R}_{H^1}(m) = \frac{\gamma_{H^1}}{2} \int_{\Omega} \|\nabla m(x)\|_2^2 \, dx,
\end{equation}
with weighting $\gamma_{H^1} > 0$,
where $\|\cdot\|_2$ denotes the Euclidean ($\ell^2$) norm of a vector in $\mathbb{R}^{d}$.
The smoothing effect introduced by $H^1$ regularization tends to blur sharp edges. 

If it is expected \textit{a priori} that there are sharp edges (i.e., voids, fibers, particles) then total variation (TV) regularization can be introduced. 
Unlike $H^1$, TV regularization is defined using a $L^1(\Omega)$ penalization on the gradient
\begin{equation}\label{eq:true_tv}
	\mathcal{R}_{TV}(m) = \gamma_{TV} \int_{\Omega} \|\nabla m(x) \|_{2}\, dx,
\end{equation}
with weighting $\gamma_{TV} > 0$.
Since the $L^1(\Omega)$ norm is sparsifying, TV regularization tends to result in piecewise constant parameter fields. 
However, $\mathcal{R}_{TV}(m)$ is not differentiable when $\nabla m=0$. 
A common approach to addressing this issue is by introducing a smoothing parameter, 
$\epsilon > 0$, to define a smoothed TV regularization term,
\begin{equation}\label{eq:epsilon_tv}
	\mathcal{R}_{TV}^\epsilon(m) = \gamma_{TV}^\epsilon \int_{\Omega} \sqrt{\|\nabla m(x)\|_2^2 + \epsilon} \, dx,
\end{equation}
with weighting $\gamma_{TV}^{\epsilon} > 0$. 
This smooths the TV functional at zero, making it differentiable and also leads to a more positive definite Hessian of $\mathcal{R}_{TV}^{\epsilon}$.
However, this will tend to produce smooth transitions between piecewise constant regions, 
where the sizes of the transition zones increase with $\epsilon$. 
Reducing $\epsilon$ brings the approximation $\mathcal{R}^{\epsilon}_{TV}$ closer to $\mathcal{R}_{TV}(m)$, 
effectively sharpening edges, but also makes the Hessian less positive-definite. 

For the Tikhonov regularization terms, one tunes only the weighting factor ($\gamma_{L^2}$ or $\gamma_{H^1}$) as a hyperparameter, whereas in the case of smoothed TV regularization, there are two hyperparameters, $\gamma_{TV}$ and $\epsilon$. In this study, $\epsilon$ is systematically reduced to sharpen the edges of the parameter field, until the problem becomes ill-conditioned again, and the optimization algorithm fails to converge. Typically, one will use a combination of regularization terms, such as combining $L^2$ with $H^1$ or TV, leading to two or three tunable hyperparameters. A common approach is to use a heuristic, such as the L-curve criterion, to determine the optimal values of the regularization parameters. Regardless of method, it is important to tune the regularization to result in a well-posed optimization problem that considers the physical understanding of the problem at hand.

\subsection{Inexact Newton-CG Method for the Regularized Inverse Problem}

The incorporation of the regularization terms into the inverse problem results in a well-posed optimization problem in which we solve with a scalable and efficient inexact Newton-CG method. To this end, we opt for a reduced-space approach, in which the minimization problem \eqref{eq:idic_minimization_problem} is solved by explicitly eliminating
the PDE constraint $R(u,m,t) = 0$ using the solution operator $u = u(m,t)$. 
This leads to an unconstrained minimization problem, 
\begin{equation}\label{eq:idic_minimization_problem_reduced}
	\min_{m \in \mathcal{M}} \widehat{\mathcal{J}}(m) := \mathcal{J}(u(m,t); I_0, I_1) = \Phi(u(m,t); I_0, I_1) + \mathcal{R}(m),
\end{equation}
since $u$ directly depends on $m$ through the solution of the PDE. To solve \eqref{eq:idic_minimization_problem_reduced}, we opt for a gradient based optimization scheme. Notably, gradient descent is not mesh-independent since it is not affine-invariant. Instead, we consider Newton's method for solving \eqref{eq:idic_minimization_problem_reduced} since its affine-invariance leads to dimension independent performance \cite{Heinkenschlosst1993}.
Starting at an initial guess, $m^0$, $m^k$ is updated at each step for $k = 1,2,...,$ by
\begin{equation}
m^{(k+1)}=m^{(k)} + \beta^{(k)} \delta m^{(k)}.
\end{equation}
where $\beta^{(k)}$ is the step size. The search direction, $\delta m^{(k)}$, is determined through the Newton step, written below in variational form,
\begin{equation}\label{eq:newton_step}
	\left\langle D_{mm} \widehat{\mathcal{J}}(m^{(k)}) \delta m^{(k)},\tilde{m} \right\rangle_\mathcal{M} = - \left\langle D_m \widehat{\mathcal{J}}(m^{(k)}),\tilde{m} \right\rangle_\mathcal{M} \qquad \forall \tilde{m} \in \mathcal{M},
\end{equation}
where $D_{m} \widehat{\mathcal{J}}$ and $D_{mm} \widehat{\mathcal{J}}$ are the first and second G\^ateaux derivatives (i.e., gradient and Hessian respectively) of the (reduced-space) cost functional. These derivatives, which are implicitly defined through the PDE solution operator, can be computed efficiently using the adjoint method; we derive their explicit forms in section \ref{subsection:derivative_derivations}.
However, explicitly constructing the full Hessian is intractable since a pair of linearized forward and adjoint PDEs must be solved to form each column of the Hessian via Hessian-vector products. The number of columns scales with the discretization dimension, not the \emph{inherent} dimension of the problem. 
Instead, we leverage the inexact Newton-CG (INCG) algorithm \cite{akccelik2006parallel, borzi2011computational}, 
which is a matrix-free Kyrlov method for the approximation solution of the Newton system \eqref{eq:newton_step}, only requiring Hessian actions, thereby avoiding explicitly forming and inverting the Hessian. Using the gradient norm for the inexact solves,
\begin{equation}
	\left\| D_m \widehat{\mathcal{J}}(m^{(k)}) + D_{mm} \widehat{\mathcal{J}}(m^{(k)}) \delta m^{(k)} \right\|_{\mathcal{M}'} \leq r_{tol} \left\| D_m \widehat{\mathcal{J}}(m^{(k)})\right\|_{\mathcal{M}'}
\end{equation}
allows one to avoid expensive computations far away from an optimum where only descent is needed. As the iterates proceed the tolerance grows tighter and better Newton iterates are produced. We use the tolerance, $r_{tol} = \text{min} \left(0.5,\sqrt{\widehat{\mathcal{J}}(m^{(k)})/\widehat{\mathcal{J}}(m^{(0)})} \right)$, as this theoretically allows for superlinear convergence in local regimes \cite{Eisenstat1994Globally}. 
Thus, with the chosen CG termination criterion, 
inexact Newton-CG ensures that only few Hessian actions are used to solve for the search direction in the global (pre-asymptotic) regime, where the descent does not benefit significantly from Hessian information. On the other hand, in the regime of local convergence, \eqref{eq:newton_step} is solved to a finer tolerance, incorporating more accurate curvature information which classically leads to superlinear convergence when the true Hessian is used. In this work we utilize the Gauss--Newton Hessian instead due to our assumed smoothness of the image field, so we settle for fast linear convergence.


\subsection{Derivation of the Adjoints and Hessian actions for $\infty$--IDIC} \label{subsection:derivative_derivations}
In this section, we derive the gradient and Hessian actions of the (reduced-space) cost functional that are required for the INCG algorithm.
To compute the gradient of $\mathcal{J}$, we employ a Lagrangian approach by introducing an adjoint variable, $p \in \mathcal{V}$ and a Lagrangian functional,
\begin{equation} \label{eq:lagrangian}
	\mathcal{L}(u,m,t,p) = \Phi(u; I_0,I_1) + \mathcal{R}(m) + r(u,m,t,p).
\end{equation}
The $m$-reduced gradient is formed by eliminating the Karush--Kuhn--Tucker (KKT) conditions associated with $p$ and $m$, i.e., we set the first variations of the Lagrangian with respect to $p$ and $u$ to zero, such that the first derivative of the cost functional with respect to $m$, $D_m \mathcal{J}$, coincides with $\partial_{m} \mathcal{L}$.
For brevity, we will omit the arguments to $\mathcal{L}$, where it is understood that they are to be evaluated at $(u,m,t)$. 
The forward equation, adjoint equation, and gradient equation are thus given by,
\begin{subequations}
\begin{alignat}{2}
	\label{eq:forward_equation}
	r(u,m,t,\tilde{p}) &= 0 
	&& \quad \forall \tilde{p} \in \mathcal{V}, \\
	\label{eq:adjoint_equation}
	\langle \partial_u r(u,m,t,p), \tilde{u} \rangle_{\mathcal{V}}
	&= - \langle \partial_u \Phi, \tilde{u} \rangle_{\mathcal{V}} 
	&& \quad \forall \tilde{u} \in \mathcal{V}, \\
	\label{eq:gradient_equation}
	\langle D_m \mathcal{J}, \tilde{m} \rangle_{\mathcal{M}}
		&= \langle \partial_m \mathcal{R}(m), \tilde{m} \rangle_{\mathcal{M}} + \langle \partial_m r(u,m,t,p), \tilde{m} \rangle_{\mathcal{M}} 
	&& \quad \forall \tilde{m} \in \mathcal{M},
\end{alignat}
\end{subequations}
respectively, where we can write out the misfit's contribution to the adjoint equation as
\begin{equation} \label{eq:misfit_gradient}
	\langle \partial_u \Phi, \tilde{u} \rangle_\mathcal{V} 
= \int_{\Omega} (I_1(x + u(x)) - I_0(x)) \nabla I_1(x + u(x)) \cdot \tilde{u}(x) \, dx  \quad \forall \tilde{u} \in \mathcal{V}.
\end{equation}
The gradient of the cost function can therefore be found by the subsequent process of first solving \eqref{eq:forward_equation} for $u$, then given $u$ solving \eqref{eq:adjoint_equation} for $p$, then given $u$ and $p$ forming the variational gradient in \eqref{eq:gradient_equation}. To derive the Hessian action, we introduce the Hessian meta-Lagrangian, 
\begin{equation} \label{eq:hessian_lagrangian}
	\mathcal{L}^H(u,m,t,p, \hat{u}, \hat{m}, \hat{p}) = \langle \partial_{u} \mathcal{L}, \hat{u} \rangle_{\mathcal{V}} 
		+ \langle \partial_{m} \mathcal{L}, \hat{m} \rangle_{\mathcal{M}} 
		+ \langle \partial_{p} \mathcal{L}, \hat{p} \rangle_{\mathcal{V}}.
\end{equation}
Analogous to the gradient formation procedure, we set $\partial_p \mathcal{L}^H = 0$, $\partial_u \mathcal{L}^H = 0$, such that $\langle D_{mm} \mathcal{J} \hat{m}, \tilde{m} \rangle_{\mathcal{M}} = \langle \partial_{m} \mathcal{L}^{H}, \tilde{m} \rangle_{\mathcal{M}}$, where $D_{mm} \mathcal{J} \hat{m}$ is the Hessian acting in a direction $\hat{m} \in \mathcal{M}$. 
These give rise to the incremental forward, incremental adjoint, and the Hessian action equations,

\begin{subequations}
\begin{alignat}{2}
	\label{eq:incremental_forward}
	\langle \partial_u r(u,m,t,\tilde{p}), \hat{u} \rangle_{\mathcal{V}} 
		&= -\langle \partial_{m} r(u,m,t,\tilde{p}), \hat{m} \rangle_{\mathcal{M}} 
	&& \quad \forall \tilde{p} \in \mathcal{V}, \\
\label{eq:incremental_adjoint}
	\langle \partial_{u}  r(u,m,t,\hat{p}), \tilde{u} \rangle_{\mathcal{V}}
		&= -\left( \langle \partial_{uu} \Phi \hat{u}, \tilde{u} \rangle_{\mathcal{V}} 
		+ \langle \partial_{mu} r(u,m,t,p) \hat{m}, \tilde{u} \rangle_{\mathcal{V}} 
		+ \langle \partial_{uu} r(u,m,t,p) \hat{u}, \tilde{u} \rangle_{\mathcal{V}}  \right)
		&& \quad \forall \tilde{u} \in \mathcal{V}, \\
\label{eq:hessian_action} 
	\langle D_{mm} \mathcal{J} \hat{m}, \tilde{m} \rangle_{\mathcal{M}}
	&= \langle \partial_{mm} \mathcal{R}(m) \hat{m}, \tilde{m} \rangle_{\mathcal{M}} 
		+ \langle \partial_{mm} r(u,m,t,\hat{p}) \hat{m}, \tilde{m} \rangle_{\mathcal{M}} \nonumber \\
		& \qquad + \langle \partial_{m} r(u,m,t,\hat{p}) \tilde{m}, \hat{p} \rangle_{\mathcal{M}} 
		+ \langle \partial_{um} r(u,m,t,p) \hat{u}, \tilde{m} \rangle_{\mathcal{M}} 
		&& \quad \forall \tilde{m} \in \mathcal{M},
\end{alignat}
\end{subequations}
respectively.

The incremental forward equation \eqref{eq:incremental_forward} is solved for $\hat{p}$, 
the incremental adjoint problem \eqref{eq:incremental_adjoint} is solved for $\hat{u}$, 
and the variables are combined with $u$ and $p$ 
to give the Hessian of the cost functional acting in the direction of $\hat{m}$ using \eqref{eq:hessian_action}.
Again, we can explicitly write out the misfit contribution to the incremental adjoint equation as 
\begin{align*}
	\langle \partial_{uu} \Phi \hat{u},  \tilde{u} \rangle_\mathcal{V} 
		&= \int_{\Omega_0} (I_1(x + u(x)) - I_0(x)) 
		\nabla^2 I_1(x + u(x)) \hat{u} \cdot  \tilde{u}  \, dx \\
		& \qquad + \int_{\Omega_0} (\nabla I_1(x + u(x)) \cdot \hat{u}) (\nabla I_1(x + u(x)) \cdot \tilde{u})  \, dx 
		\quad \forall \hat{u}, \tilde{u} \in \mathcal{V}.
\end{align*}
The second spatial derivative of an image, specifically $\nabla^2 I_1(x + u(x)) $, may not exist given our assumed $\mathcal{C}^1$ smoothness.
To avoid more restrictive assumptions of $\mathcal{C}^2$ smoothness, we employ a Gauss--Newton approximation of the Hessian \cite{Heinkenschlosst1993},
in which we drop the $(I_0(x) - I_1(x + u(x))) \nabla^2 I_i(x + u(x)) \hat{u} \cdot  \tilde{u}$ term, leaving only the terms involving $\nabla I_1$. 
That is, we take the approximation 
\begin{equation}\label{eq:gauss_newton_hessian}
	\langle \partial_{uu} \Phi \hat{u},  \tilde{u} \rangle_\mathcal{V} 
		\approx \int_{\Omega_0} (\nabla I_1(x + u(x)) \cdot \hat{u}) (\nabla I_1(x + u(x)) \cdot \tilde{u})  \, dx 
		\quad \forall \hat{u}, \tilde{u} \in \mathcal{V}.
\end{equation}
This approximation is accurate when the data misfit $I_1(x + u(x)) - I_0(x)$ is small, 
but since $I_1 \not\to I_0$ due to the observational noise, the Gauss--Newton Hessian does not converge to the true Hessian. We do not expect the typical superlinear convergence of inexact Newton methods and instead fast linear convergence. 

\begin{remark}
	It is possible that the true Hessian can be computed if the second spatial derivative of the image is well-defined, but this is not the case for the images considered in this work as they are piece-wise constant speckle patterns where the second derivative is undefined at the edges of the speckles.
\end{remark}

\subsection{Derivations of the gradient and Hessian action for Regularization Terms}

In this section we derive the gradient and Hessian actions of the regularization terms. 
Beginning with the $L^2$ Tikohnov regularization, we have the following gradient and Hessian action:
\begin{subequations}
\begin{alignat}{2}
\label{eq:grad_l2}
	\langle D_{m} \mathcal{R}_{L^2}(m), \tilde{m} \rangle_{\mathcal{M}} 
		&= \gamma_{L^2} \int_{\Omega} m \tilde{m} \, dx
		&& \quad \forall \tilde{m} \in \mathcal{M}, \\
\label{eq:hessian_action_l2}
	\langle D^2_{mm} \mathcal{R}_{L^2}(m) \hat{m}, \tilde{m} \rangle_{\mathcal{M}}
	&= \gamma_{L^2} \int_{\Omega} \hat{m} \tilde{m} \, dx
		&& \quad \forall \tilde{m}, \hat{m} \in \mathcal{M}.
\end{alignat}
\end{subequations}
For the $H^1$ regularization, we have the following gradient and Hessian action:
\begin{subequations}
\begin{alignat}{2}
\label{eq:grad_h1}
	\langle D_{m} \mathcal{R}_{H^1}(m), \tilde{m} \rangle_{\mathcal{M}} 
		&= \gamma_{H^1} \int_{\Omega} \nabla m \cdot \nabla \tilde{m} \, dx
		&& \quad \forall \tilde{m} \in \mathcal{M}, \\
\label{eq:hessian_action_h1}
	\langle D^2_{mm} \mathcal{R}_{H^1}(m) \hat{m}, \tilde{m} \rangle_{\mathcal{M}}
	&= \gamma_{H^1} \int_{\Omega} \nabla \hat{m} \cdot \nabla \tilde{m} \, dx
		&& \quad \forall \tilde{m}, \hat{m} \in \mathcal{M}.
\end{alignat}
\end{subequations}
For the TV regularization we have the following gradient and Hessian action, which we refer to as the primal formulation of TV:
\begin{subequations}\label{eq:tv_derivatives}
\begin{alignat}{2}
\label{eq:grad_tv}
	\langle D_{m} \mathcal{R}_{TV}^\epsilon(m), \tilde{m} \rangle_{\mathcal{M}} 
	 & = \gamma_{TV}^\epsilon \int_{\Omega} \frac{\nabla m}{\sqrt{\|\nabla m\|_2^2 + \epsilon}} \cdot \nabla \tilde{m} \, dx
	  && \quad \forall \tilde{m} \in \mathcal{M}. \\
\label{eq:hessian_action_tv}
	\langle D^2_{mm}  \mathcal{R}_{TV}^\epsilon(m) \hat{m}, \tilde{m} \rangle_{\mathcal{M}}
		&= \gamma_{TV}^\epsilon \int_{\Omega} \frac{1}{\sqrt{\|\nabla m\|_2^2 + \epsilon}} \left[ \left(I - \frac{\nabla m \otimes \nabla m}{ \nabla m \nabla m + \epsilon}\right) \nabla \hat{m} \right] \cdot \tilde{m} \, dx
		&& \quad \forall \tilde{m}, \hat{m} \in \mathcal{M},
\end{alignat}
\end{subequations}
where $\|\cdot\|_2$ denotes the usual Euclidean norm for vectors in $\mathbb{R}^d$. It has been shown that anisotropy in the $(I - \frac{\nabla m \otimes \nabla m}{ \nabla m \nabla m + \epsilon})$ term causes convergence issues; as $\epsilon$ is reduced the radius of local convergence also is reduced and the method may not converge globally \cite{Chan1999}. Additionally, the $\epsilon$ parameter is mesh dependent, so one must retune the hyperparameter as the mesh is refined. A first attempt to mitigate this is to scale $\epsilon$ with the mesh size \cite{Chan1999_lag}. However, this method proved unsuccessful in our experiments and mesh-independence was not achieved, possibly due to the nonlinear nature of the $\infty$--IDIC formulation (and nonlinear PDE used in demonstration). To mitigate the mesh-depenence, we instead employ a primal-dual approach which empirically leads to a mesh-independent algorithm.

\subsubsection{Primal-Dual TV}


The issues with convergence of Newton's method for the primal $\mathcal{R}^\epsilon_{TV}$ are due to the highly ill-conditioned Hessian term \eqref{eq:hessian_action_tv}. The ill-conditioning of the Hessian is due to taking an additional derivative of the $\frac{\nabla m}{\sqrt{\|\nabla m\|^2_m + \epsilon}}$ term in the gradient \eqref{eq:grad_tv}. In \cite{Chan1999}, the authors proposed a solution to this problem via an equivalent \emph{primal-dual} formulation of the TV-regularized optimization problem. In this formulation an additional dual variable $w= \frac{\nabla m}{\sqrt{\|\nabla m\|^2_m + \epsilon}}$ is introduced leading to the primal-dual interpretation of TV:
\begin{equation}\label{eq:abstract_primal_dual_tv}
	\mathcal{R}^\epsilon_{TV}(m) = \gamma_{TV}\sup_{w \in \mathcal{W}}\left\{\int_\Omega m(x) \nabla \cdot w(x) dx: \|w(x)\|_2 \leq 1 \enskip \forall x \in \Omega\right\},
\end{equation}
where $\mathcal{W}$ is an appropriately chosen function space for the dual variable. For example when approximating $\mathcal{M}$ by continuous Galerkin finite elements with order $r$, $\mathcal{C}\mathcal{G}^r(\Omega)$  the space of ${r-1}^{th}$ order discontinuous Galerkin finite elemets, $\mathcal{DG}_r(\Omega)^d$, is a suitable choice for $\mathcal{W}$, \cite{arnold2010finite,herrmann2019discrete}.
In the primal-dual TV, the first order optimality condition \eqref{eq:gradient_equation} for $m$,  
\begin{equation}
	\langle D_m \mathcal{J}, \tilde{m} \rangle_{\mathcal{M}} =0\qquad\forall \tilde{m} \in \mathcal{M}
\end{equation}
can be rewritten as a system for the primal-dual pair $(m,w)$:
\begin{subequations}
\begin{align}\label{eq:primaldual_optimality}
\gamma_{TV}^\epsilon \int_{\Omega} w \cdot \nabla \tilde{m} \, dx + \langle \partial_m r(u,m,t,p), \tilde{m} \rangle_{\mathcal{M}} &= 0\quad \forall \tilde{m} \in \mathcal{M}, \\
\int_\Omega\left(\sqrt{\|\nabla m\|_2^2 + \epsilon}\right) w \tilde{w}dx - \int_\Omega \nabla m \tilde{w}dx &= 0 \quad\forall \tilde{w} \in \mathcal{W},
\end{align}
\end{subequations}
The subsequent linearization to a Newton system for the pair $(m,w)$ leads to a better-conditioned optimization problem and empirically superior convergence properties \cite{Chan1999}. At every Newton iteration, we compute updates for both $m$ and $w$, as follows:
\begin{align*}
m^{(k+1)} &= m^{(k)} +\alpha_m^{(k)}\delta m^{(k)},\\
	w^{(k+1)} &= w^{(k)} +\alpha_w^{(k)}\delta w^{(k)}.
\end{align*}
We first solve for $\delta m^{(k)}$ via the Newton system \eqref{eq:newton_step}, where the Hessian vector product is defined as in \eqref{eq:hessian_action}, with
\begin{subequations}
\begin{align}
	\langle \partial_{mm} \mathcal{R}_{TV}^\epsilon(m^{(k)}) \delta m^{(k)}, \tilde{m} \rangle_{\mathcal{M}}  &= \gamma_{TV}^\epsilon \int_\Omega \frac{1}{\sqrt{\|\nabla m^{(k)}\|_2^2 + \epsilon}} \left[ \left( I - A(m^{(k)},w^{(k)})\right) \nabla \tilde{m}\right] \cdot \nabla \delta m^{(k)}\, dx,\\
	A(m,w) &= \frac{1}{2} w \otimes \frac{\nabla m}{\sqrt{\|\nabla m\|_2^2 + \epsilon}} + \frac{1}{2} \frac{\nabla m}{\|\sqrt{|\nabla m\|_2^2 + \epsilon}} \otimes w.
\end{align}
\end{subequations}
Once $\delta m^{(k)}$ is computed, $m^{(k+1)}$ is found via a line search to determine an appropriate choice of $\alpha_m^{(k)}$. Next $\delta w^{(k)}$ is found from the following relation:
\begin{equation}\label{eq:what_pdtv}
	\langle \delta w^{(k)},\tilde{w}\rangle_{\mathcal{W}} = \int_\Omega \frac{1}{\sqrt{\|\nabla m^{(k+1)}\|_2^2 + \epsilon}}\left(I - A(m^{(k+1)},w^{(k)}) \right)\nabla \tilde{w} - {w}^{(k)}\tilde{w} + \frac{\nabla m^{(k+1)}}{\sqrt{\|\nabla m^{(k+1)}\|_2^2 + \epsilon}}\tilde{w}\, dx \quad \forall \tilde{w} \in \mathcal{W}.
\end{equation}
The dual variable $w^{(k+1)}$ is then updated via a line search for $\alpha_w^{(k)}$, which is chosen such that the condition $\|w(x)|\_2\leq 1$ is met for all $x\in \Omega$, as in \eqref{eq:abstract_primal_dual_tv}. Similar to \cite{Chan1999} we observe significantly improved performance over the primal version of TV. In particular, our choice of $\epsilon$ is mesh-independent, and we observe consistent and improved convergence across all meshes. With the primal version of TV we were only able to observe convergence on small meshes, after significant hyperparameter tuning to find an appropriate choice of $\epsilon$.


\section{Numerical Results}

\subsection{Introducing the PDEs}

We evaluate $\infty$--IDIC's efficacy via synthetic data from two PDE models for elastic deformation problems: linear elasticity and hyperelasticity. The parameter fields are complex, spatially varying examples which represent materials with voids, stiff particles, and grain structures. We present the PDEs in the following section, but leave details of the data generation (i.e., image creation and noise corruption) to \ref{appendix:a}.

\subsubsection{Linear Elasticity}

We consider unit square physical domains $\Omega = (0,1)^2$, with spatial coordinates $x\in \Omega$. The strong form of the linear elasticity PDE over $\Omega$ is given by 
\begin{equation} \label{eq:linear_elasticity_strong_form}
    \begin{aligned}
        \nabla \cdot (\sigma(u)) &= 0 & \quad &\text{in} \quad \Omega, \\
        u &= 0 & \quad &\text{on} \quad \Gamma_{L}, \\
        \sigma(u) \cdot n &= t & \quad &\text{on} \quad \Gamma_{R}, \\
        \sigma(u) \cdot n &= 0 & \quad &\text{on} \quad \Gamma_{T} \cup \Gamma_{B}.
    \end{aligned}
\end{equation}
We prescribe a fixed displacement condition on the left boundary 
$\Gamma_{L}$, 
a traction $t$ on the right boundary $\Gamma_{R}$,
and traction-free boundary conditions on top and bottom boundaries, $\Gamma_{T}$ and $\Gamma_{B}$ (see Figure \ref{fig:deformation_schematic}).
Here, the stress tensor is given by 
$\sigma(u) = \lambda \nabla \cdot u I + 2 \mu \varepsilon(u)$, where
$\varepsilon(u) = \frac{1}{2}(\nabla u + \nabla u^T)$,
and $\lambda$ and $\mu$ are the Lamé parameters.
We write the traction in terms of its normal and shear components,
\begin{equation} \label{eq:traction_condition}
    t = t_{\text{normal}} e_1 + t_{\text{shear}} e_2.
\end{equation}

The Lamé parameters are related to the Young's modulus, $E$, and Poisson's ratio, $\nu$, by
\begin{equation} \label{eq:lame_parameters}
    \lambda = \frac{E \nu}{(1 + \nu)(1 - 2\nu)}, \quad \mu = \frac{E}{2(1 + \nu)}.
\end{equation}
We fix the Poisson's ratio as $\nu = 0.35$ and invert for the Young's modulus, $E(x)$.
To prevent the modulus from becoming negative, we parametrize it using the exponential function, i.e., $E(x) = e^{m(x)}$, and invert for the log-Young's modulus $m$. Using the notation $H^1_{L}(\Omega) = \{ u \in H^1(\Omega) : u|_{\Gamma_L} = 0 \}$, 
the variational form of the PDE \eqref{eq:forward_pde_weak} is: find \( u \in H^1_L(\Omega) \) such that,
\begin{equation} \label{eq:linear_elasticity_variational_formulation}
   r_{\mathrm{LE}}(u,m,t,v) := \int_{\Omega} \sigma(u, m) : \varepsilon(v) \, dx - \int_{\Gamma_L} t \cdot v \, ds = 0 \qquad \forall v \in H^1_L(\Omega).
\end{equation}

\subsubsection{Neo-Hookean Hyperelasticity}
As a second PDE problem we consider a hyperelastic material model, 
which is commonly used to model rubber-like materials, biological tissues, and other materials that undergo large deformations. This problem setting is of interest for elastography \cite{goenezen2011solution, mei2016estimating}. The stress-strain relationship is nonlinear in nature, making it a more difficult numerical test case. Here we use a material point description, unlike linear elasticity where we use spatial points. Let $X\in \Omega$ be the material point where $X$ is related to the spatial point, $x$, through the displacement, $u(X) = x -X$. In hyperelasticity, the internal forces that develop in the material are derived from a strain energy function, $W = W(X,\mathbf{C}(X))$, where $\mathbf{C}(X)=\mathbf{F}(X)^T \mathbf{F}(X)$ is the right Cauchy--Green deformation tensor and $\mathbf{F}(X)=I+\nabla u(X)$ is the deformation gradient. 
For the neo-Hookean model, the strain energy function is given by 

\begin{equation} \label{eq:neo_hookean}
    W (X,\mathbf{C}(X)) = \frac{\mu(X)}{2} (\text{tr}(\mathbf{C}(X)) - 3) - \mu(X) \ln J(X) + \frac{\lambda(X)}{2} \ln^2 J(X),
\end{equation}
where $\mu$ and $\lambda$ are the Lam\'e parameters, $\text{tr}(\mathbf{C})$ is the trace of $\mathbf{C}$, and $J = \det(F)$ is the determinant of the deformation gradient. 
The strong form of hypelasticity PDE is then given by,
\begin{equation} \label{eq:hyperelasticity_strong_form}
    \begin{aligned}
        \nabla \cdot (\mathbf{F}(u)\mathbf{S}(u)) &= 0 & \quad &\text{in} \quad \Omega, \\
        u &= 0 & \quad &\text{on} \quad \Gamma_{L}, \\
        \mathbf{F}(u)\mathbf{S}(u) \cdot n &= t & \quad &\text{on} \quad \Gamma_{R}, \\
        \mathbf{F}(u)\mathbf{S}(u) \cdot n &= 0 & \quad &\text{on} \quad \Gamma_{T} \cup \Gamma_{B}.
    \end{aligned}
\end{equation}
where $\mathbf{S}(X) = 2\frac{\partial W}{\partial \mathbf{C}(X)}$ is the second Piola--Kirchhoff stress tensor,
and the applied traction $t$ is the same as in Equation \eqref{eq:traction_condition}.

Similar to linear elasticity, the Young's modulus defined using the parametrization $E(x) = e^{m(x)}$ to ensure positivity and fix the Poisson’s ratio as $\nu$ = 0.35. 
Thus, the variational form of the hyperelasticity PDE is: find \( u \in H^1_L(\Omega) \) such that for all \( v \in H^1_L(\Omega) \),
\begin{equation} \label{eq:hyperelasticity_variational_formulation}
    r_{\text{HE}}(u,m,t,v) = \int_{\Omega} \mathbf{S}(u, m) : \varepsilon(v) \, dx - \int_{\Gamma_R} t \cdot v \, ds = 0.
\end{equation}

We use the finite element method to solve both PDE problems. In particular, we employ piecewise linear Lagrangian triangular elements to form finite element spaces $\mathcal{M}^h \approx \mathcal{M}$ and $\mathcal{V}^h \approx \mathcal{V}$. For primal-dual total variation, We discretize the dual variable, $w$, using 0th order discontinuous Galerkin to form a finite element space $\mathcal{W}^h \approx \mathcal{W}$. We utilize an $N \times N$ mesh, where $N$ is the number of elements per side. We handle all finite element calculations in FEniCS \cite{Fenics2015} and use hIPPYlib to assist with the adjoint computations \cite{Villa2021}. 


\subsection{Spatially Varying Inversions of Heterogeneous Fields}

\begin{figure}[H]
    \centering
    \includegraphics[width=.78\textwidth]{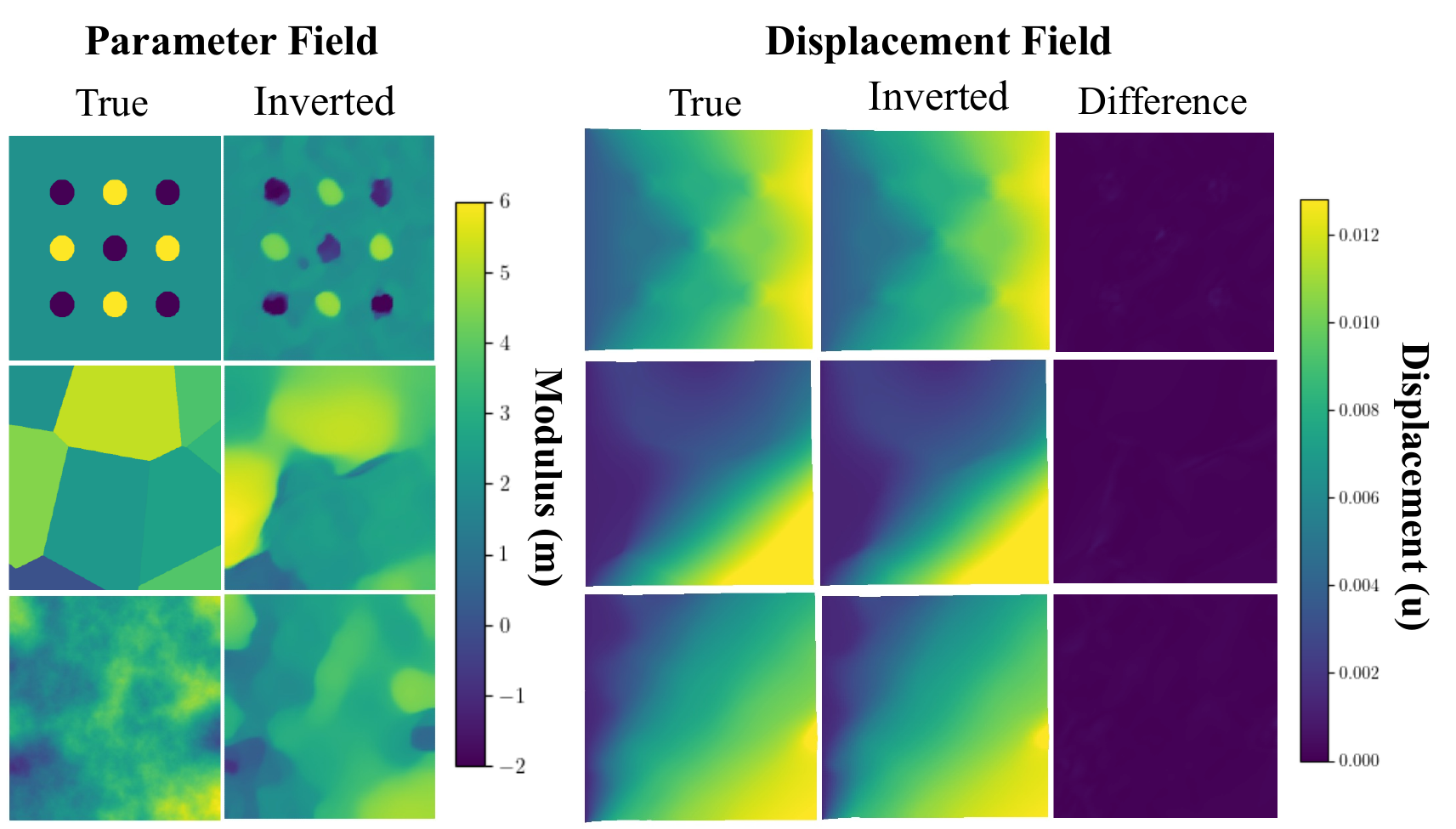}
    \includegraphics[width=.20\textwidth]{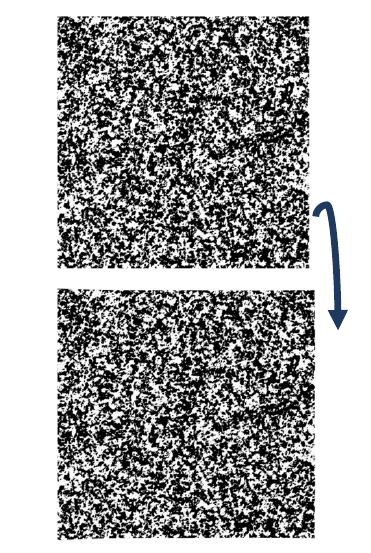}
    \caption{Inverse results for Linear Elasticity for varying features. Each inversion is based on a simple tensile experiment with $\sim$0.2\% strain. 10\% and 5\% noise were applied to the image brightness and force measurements, respectively. The inverse problem is solved using the primal-dual total variation formulation ($\gamma_{L^2} = 5\times10^{-8}$, $\gamma_{TV} = 1\times10^{-6}$) with a mesh size of $100 \times 100$, and the speckle correlation length is 0.01. The initial guess was uniformly $m(x)=2$. The synthetic images are shown on the right.}
    \label{fig:inverse_results_linear}
\end{figure}

\begin{figure}[H]
    \centering
    \includegraphics[width=.78\textwidth]{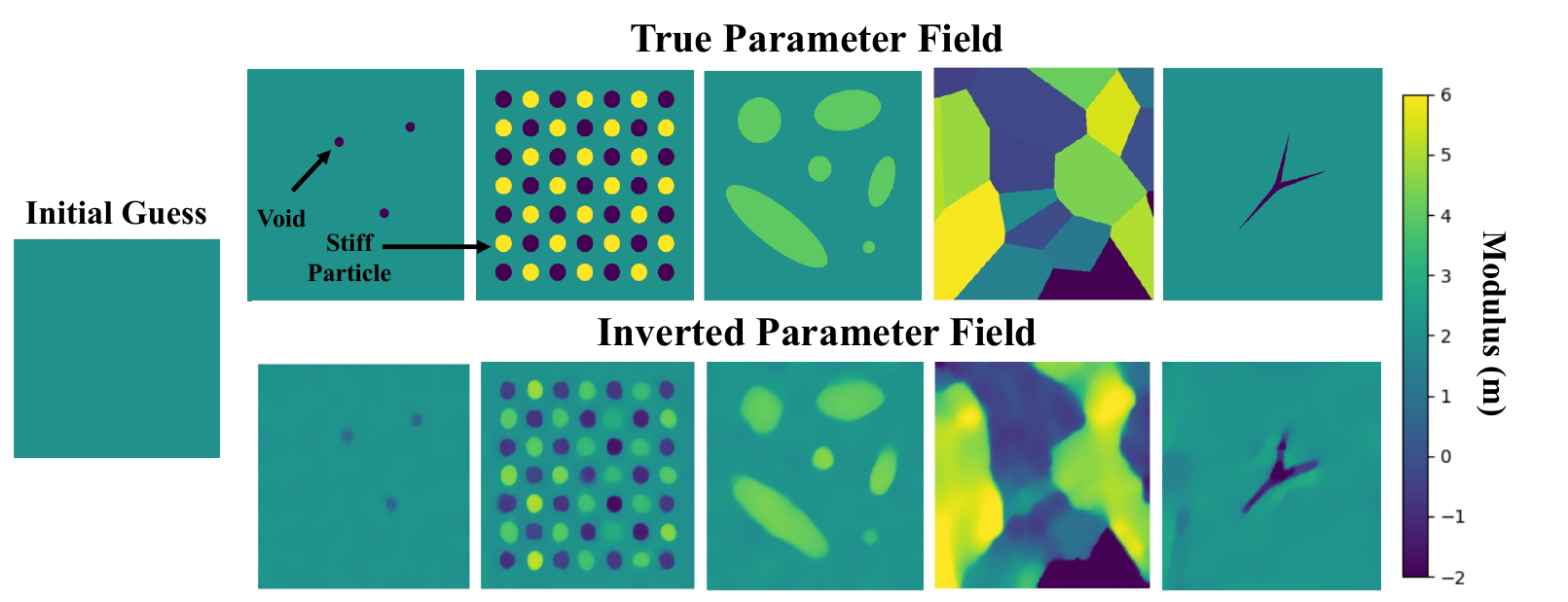}
    \includegraphics[width=.20\textwidth]{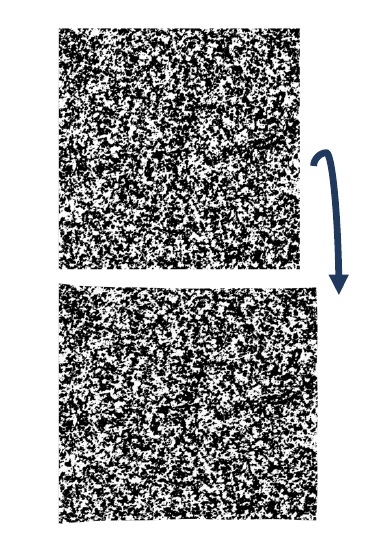}
    \caption{Inverse results for Hyperelasticity for varying features. Each inversion is based on a simple tensile experiment with $\sim$5\% strain. 10\% and 5\% noise were applied to the image brightness and force measurements, respectively. The inverse problem is solved using the primal-dual total variation formulation ($\gamma_{L^2} = 5\times10^{-6}$, $\gamma_{TV} = 7.5 \times10^{-4}$) with a mesh size of $100 \times 100$, and the speckle correlation length is 0.01.}
    \label{fig:inverse_results_hyperelasticity}
\end{figure}

We begin our numerical experiments by presenting three examples of linear elasticity generated with small strains ($\sim$0.5\%) in Figure \ref{fig:inverse_results_linear}. Even at small strains, the inversions capture the spatial distribution well. Notably, the inversions tend to prefer piece wise constant features due to the choice of TV regularization, which is particularly useful when sharp features (i.e., fibers, grains, etc.) are present, as seen in the first two examples. Additionally, we present inversions for hyperelasticity, which could pose greater difficulty due to higher deformations ($\sim$3\% strain). In Figure \ref{fig:inverse_results_hyperelasticity}, we show five inverse problems ranging from small voids and grain boundaries to cracks. It is suspected that the increase in strain led to improved inversions. Notably, we successfully find features of varying geometries and sizes, as demonstrated by complex blobs, Voronoi tessellations, and a thin, branching crack. In the first case, we detected tiny void regions ($m(x) = -2$) with a radius 1/50th of the domain length; Existing FEMU studies showed inversions for two features with a radius 1/5th of the domain length \cite{goenezen2011solution,mei2018comparative}. However, in all of our inversions, the modulus values do not perfectly match, due to the inverse problem's ill-posedness.

Ill-posedness likely arises from multiple sources including the non-uniqueness of modulus fields that produce displacement fields satisfying the image data misfit. As seen in Figure \ref{fig:inverse_results_linear}, the inverted displacement fields almost replicate the true displacement fields despite the modulus field not perfectly matching. There is a non-unique relationship between $m$ and $u$ that satisfies the data, resulting in a non-unique solution to the inverse problem. This effect may also be influenced by the initial guess ($m(x)=2$), which is two orders of magnitude different from the true solution's features. The synthetic image and force data were corrupted with noise (see \ref{appendix:a}) and noise significantly impacts the inversions due to the problem's ill-conditioning. A key aspect of the ill-posedness is the inherent dissipative nature of elasticity due to the damping of high frequency modes of $m$ by the forward physics \cite{vogel2002computational,Tarantola2005}. Unaddressed, these high frequency modes would potentially lead to unstable growth which further necessitates the need for regularization. While we incorporate regularization to reduce the ill-posedness, this ultimately introduces a modeling bias into the inversion. For example, in the case of Gaussian random field (Figure \ref{fig:inverse_results_linear}) with smooth features, the TV regularization tends to sharpen edges.

\subsection{Regularization Mitigates the Ill-Posedness}

\begin{figure}[H]
    \centering
    \includegraphics[width=.60\textwidth]{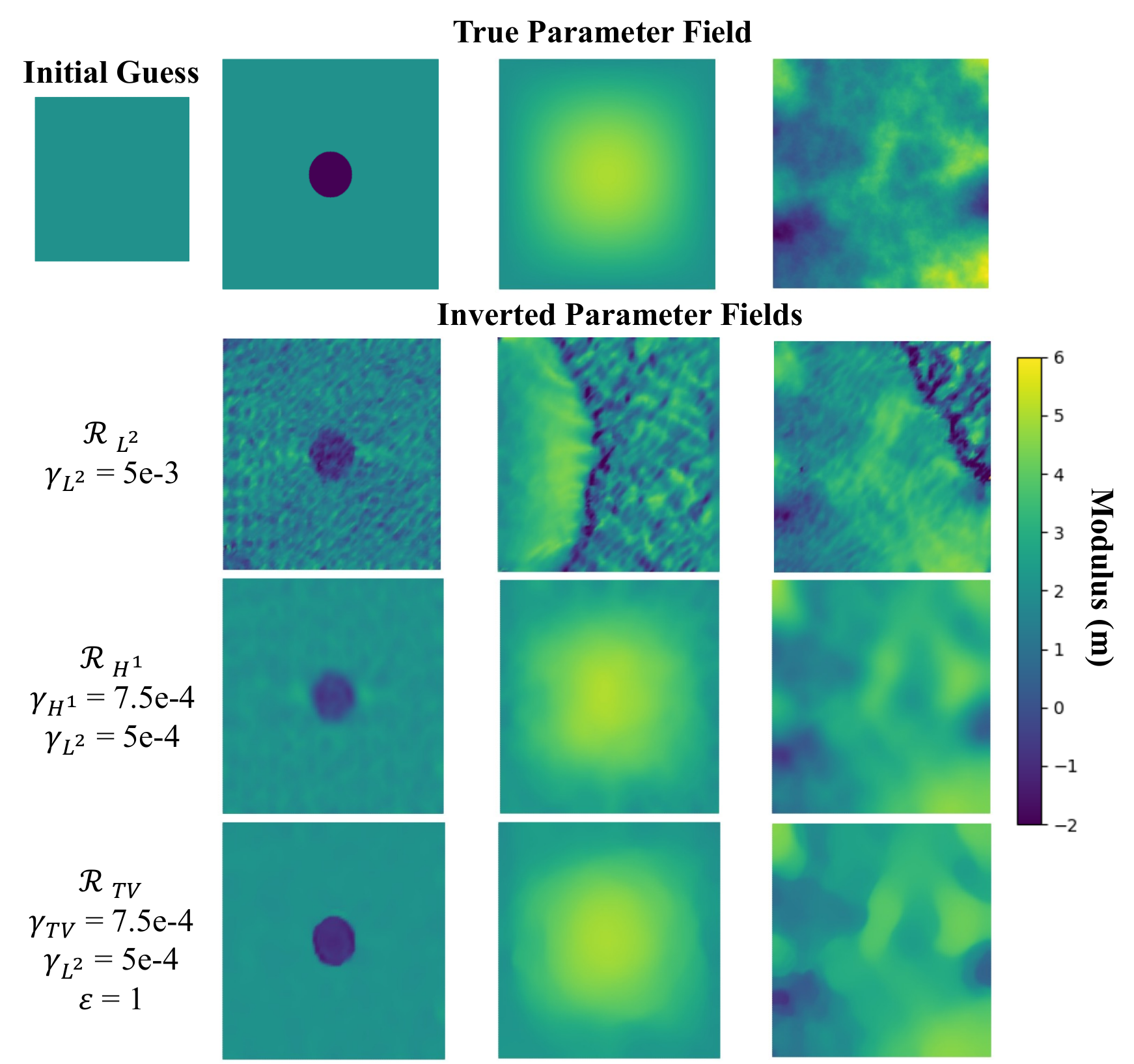}
    \caption{Inverse results for Hyperelasticity for varying regularization. Each inversion is based on a simple tensile experiment with $\sim$5\% strain. 10\% and 5\% noise were applied to the image brightness and force measurements, respectively. The inverse problem is solved using a mesh size of $100 \times 100$ and a speckle correlation length of 0.01.}
    \label{fig:regularization_grid}
\end{figure}

In this section, we study how the choice of regularization affects the inferred parameter fields. In particular, we solve the inverse problem using $L^2$, $H^1$ and TV regularizations for 
three distinct $m_{\text{true}}$ scenarios: a void inclusion, a bump function, and a Gaussian random field. 
The inverted parameter fields are shown in Figure~\ref{fig:regularization_grid}. These results demonstrate that $L^2$ regularization is insufficient and introduces artifacts, despite considerable tuning efforts for the weighting parameter, yielding noisy and inaccurate inversions incapable of identifying material features. Although $L^2$ penalizes deviations from a nominal value, and in principle addresses the ill-posedness of the inverse problem, it does not adequately capture the structure of the material fields. In contrast, both $H^1$ and TV regularization successfully identify features, with varying degrees of smoothing. TV regularization excels in preserving sharp interfaces, as seen in the void inclusion scenario, while $H^1$ regularization produces smoother inversions without favoring sharp interfaces. Notably, $H^1$ regularization is better suited than TV for the Gaussian random field, highlighting the importance of adapting regularization based on the nature of the problem. These observations underscore the significance of selecting an appropriate regularization since the choice determines the preferred modulus field characteristics when the noisy data lone is insufficient to recover the true modulus. 

\subsection{Mesh-independent Performance}

\begin{figure}[H]
    \centering
    \includegraphics[width=.75\textwidth]{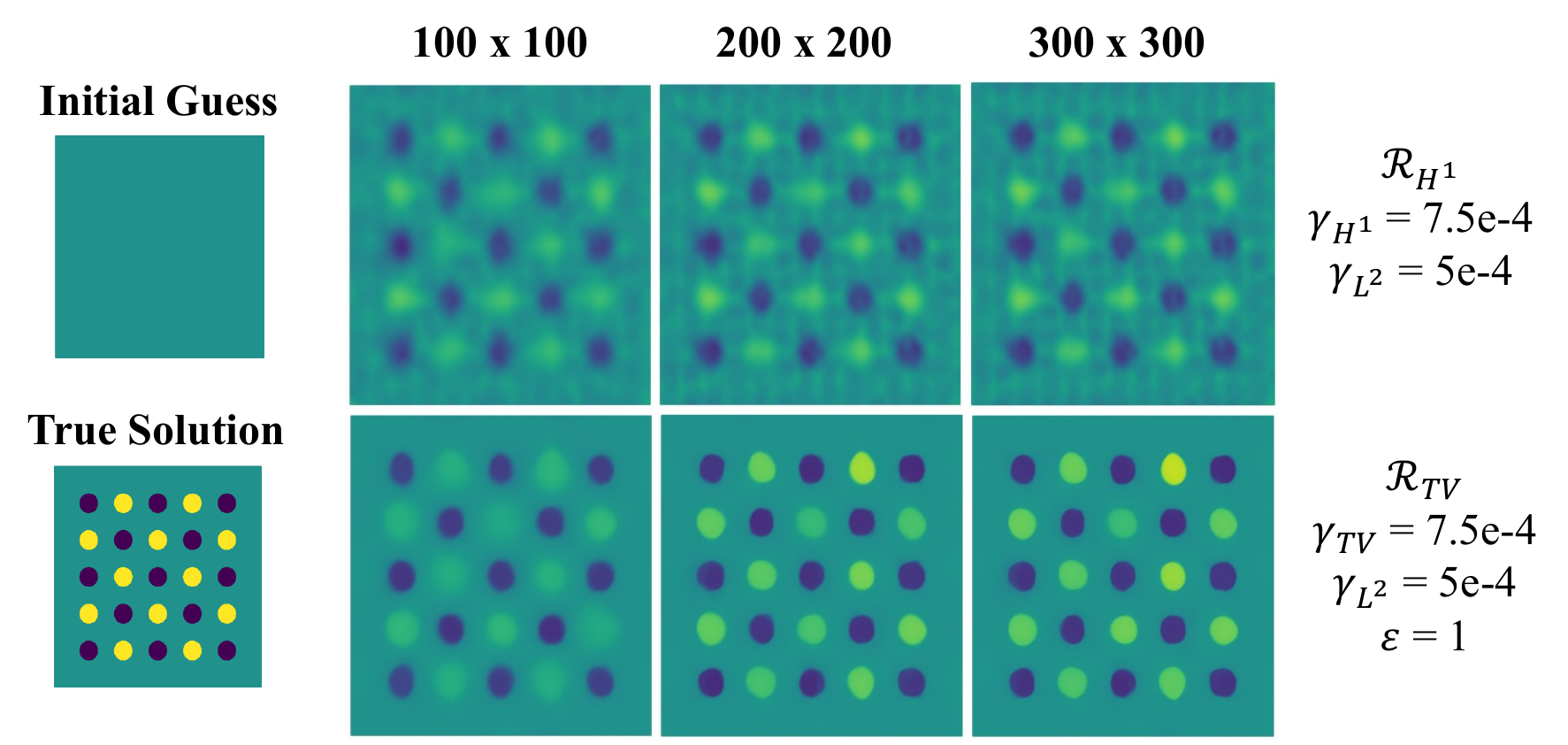}
    \caption{Performance of $\infty$--IDIC with $H^1$ and TV regularization for the sample problem with varying mesh resolution. Each inversion is based on a simple tensile experiment with $\sim$5\% strain. 10\% and 5\% noise were applied to the image brightness and force measurements, respectively. The speckle correlation length is 0.01.}
    \label{fig:mesh_size_examples}
\end{figure}

The resolution of the mesh must be chosen to sufficiently be able to resolve the physics (e.g., feature size, such as for the thin, branched crack in Figure \ref{fig:inverse_results_hyperelasticity}). We show inversions for 40,401, 160,801, and 361,201 parameter values (degrees of freedom) relating to $100 \times 100$, $200 \times 200$ and $300 \times 300$ meshes, respectively. Previous IDIC and FEMU studies handled parameter fields in low dimensions such as two \cite{Mathieu2015EstimationIntegrated-DIC}, three \cite{Neggers2019SimultaneousIdentification} and eight \cite{Gras2015} scalar values. Increasing the number of degrees of freedom leads to smaller features being captured. Additionally, fine mesh resolution is important for representing the image speckle pattern. As we refine the mesh, the speckle pattern is  better represented which leads to a better inversion as seen in Figure \ref{fig:mesh_size_examples}. The synthetic images were generated on a $500 \times 500$ mesh and then corrupted with noise; Since, the $\infty$--IDIC was conducted on coarser meshes there is some information loss in the finite element representation of the images. Refining the mesh is necessary to minimize this effect, but it comes at additional computational costs. 

\begin{figure}[H]
    \centering
    \includegraphics[width=.75\textwidth]{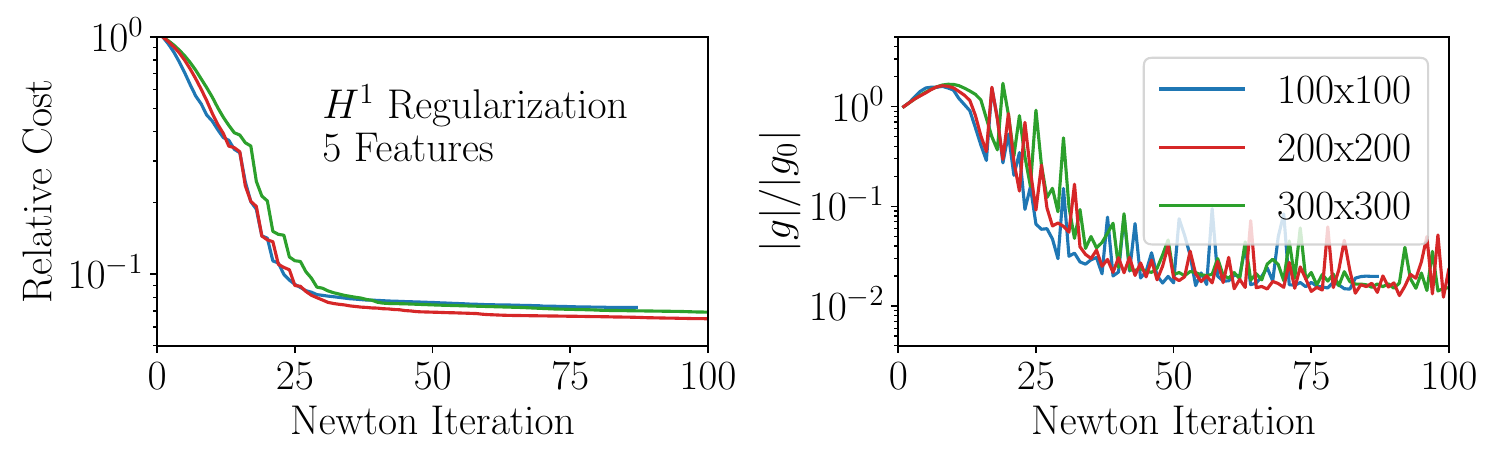}
    \caption{Mesh-independent performance of $\infty$--IDIC for $H^1$ regularization on the parameter where the mesh increased from $100 \times 100$, $200 \times 200$ to $300 \times 300$. This is a tensile experiment using a $5 \times 5$ features on hyperelastic medium with $\sim$5\% strain. 10\% and 5\% noise were applied to the image brightness and force measurements, respectively. The speckle correlation length is 0.01.}
    \label{fig:H1_mesh_independence}
\end{figure}

A key benefit of formulating the problem in infinite dimensions is mesh-independent optimization behavior, leading to a scalable inversion algorithm. Refining the mesh only increases the cost per PDE solve and not the number of optimization iterations. Figure \ref{fig:H1_mesh_independence} demonstrates the mesh-independent performance of $\infty$--IDIC for $H^1$ regularization on the parameter where the mesh increased from $100 \times 100$, $200 \times 200$ to $300 \times 300$. We show the relative cost (cost at each Newton iteration normalized by the initial cost) and the gradient norm. After 30 iterations, the relative cost plateaus independent of mesh size. The gradient norm exhibits similar convergence behavior, but the refined mesh leads to more jumps in the gradient values. 

\begin{figure}[H]
    \centering
    \includegraphics[width=.75\textwidth]{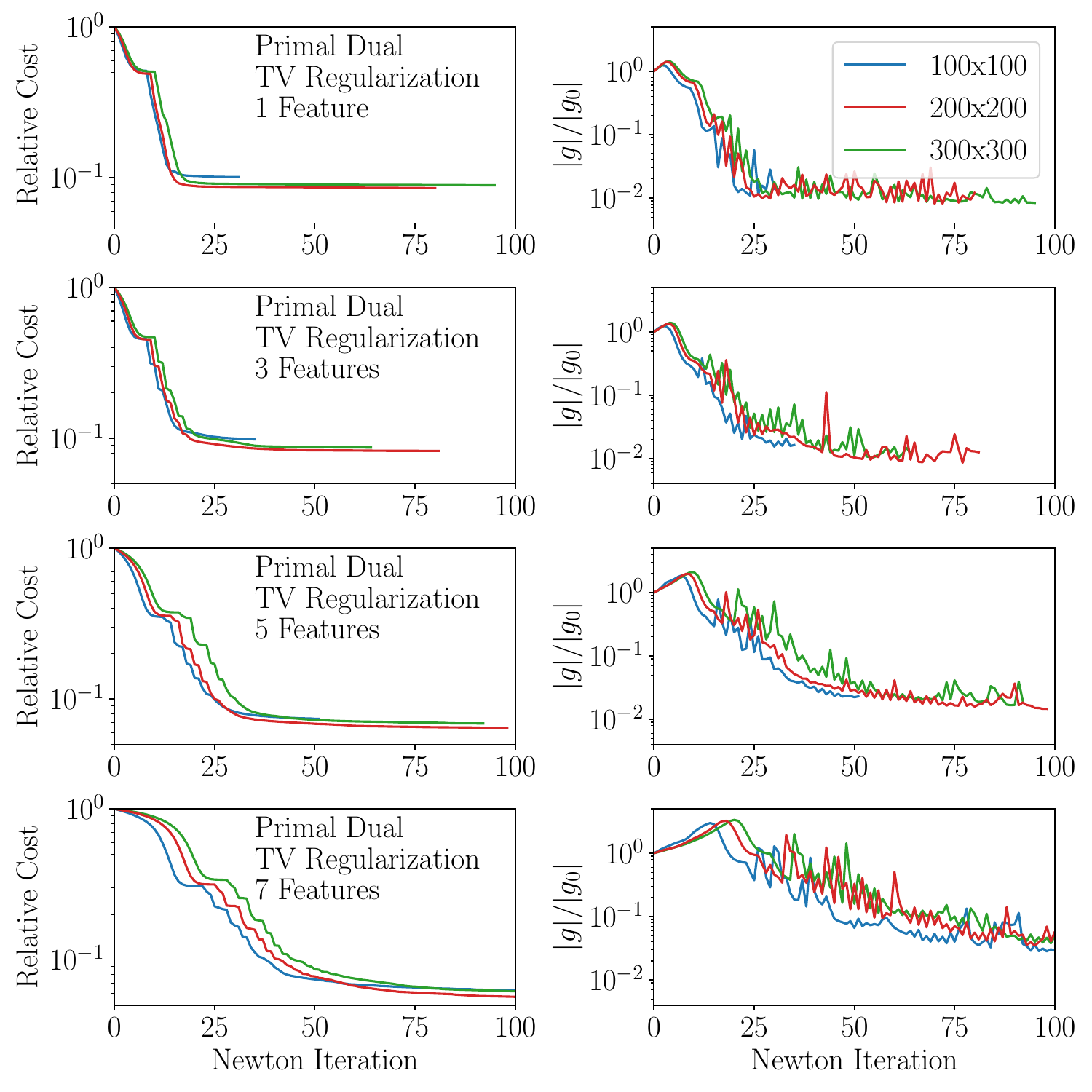}
    \caption{Mesh-independent performance of $\infty$--IDIC for Hyperelasticity when using the primal-dual total variation formulation ($\gamma_{L^2} = 5\times10^{-6}$, $\gamma_{TV} = 7.5 \times10^{-4}$) where the mesh increased from $100 \times 100$, $200 \times 200$ to $300 \times 300$.}
    \label{fig:TV_mesh_independence}
\end{figure}

We demonstrate mesh-independent performance using the primal-dual TV formulation in Figure \ref{fig:TV_mesh_independence}. It is clear from the results that $\infty$--IDIC performs independently of the mesh size. The relative cost reduction and normalized gradient norm demonstrate similar convergence behavior among the varying mesh sizes. However, varying the number of features in $m_{\text{true}}$ results in different convergence behavior, as the complexity of the problem increases. With more features, the problem requires more iterations to converge because each experiment starts from the same initial guess of $m(x) = 2$, and this guess is further from the true solution when more features are present. Additionally, the finer mesh better represents the image speckle, an effect that may be more pronounced for more complex problems with smaller features (e.g., $7 \times 7$). This effect may also be influenced by the regularization weighting. The tuning of the regularization parameter likely depends on the specific problem, and all examples in Figure \ref{fig:TV_mesh_independence} use the same regularization weighting. Ultimately, each experiment converges in a mesh-independent manner leading to a scalable algorithm that can be refined to represent the image speckle and resolve finer material features. 

\begin{table}[H]
    \centering
    \begin{tabular}{|c|c|}
        \hline
        Mesh Size & Wall Clock Time (mins)  \\
        \hline
        $100 \times 100$ & 2.20 \\
        $200 \times 200$ & 7.30 \\
        $300 \times 300$ & 35.26 \\
        \hline
    \end{tabular}
    \caption{Timing results for the three different mesh sizes. $\infty$--IDIC was run on an M2 Macbook Pro with 16 Gb of memory in serial.}
    \label{tab:timing_results}
\end{table}

Increasing the mesh increases the computation time due to the PDE solves being more expensive. Table \ref{tab:timing_results} shows the wall clock time for the three different mesh sizes. While the optimization behavior remains consistent, the computational cost increases significantly with the mesh size. The cost of the PDE solve and calculating $I(x+u(x))$ and $\nabla I(x+u(x))$ are the primary drivers of the computational cost. Notably, Pan et al. show that common DIC software (Vic-2D) took from 7 to 50 minutes to solve for the displacement field of a tensile coupon \cite{pan2011fast}. $\infty$--IDIC is competitive with these times regardless of the solving a more complex inverse problem.

\subsection{$\infty$--IDIC Robustness}

In this section, we investigate the robustness of $\infty$--IDIC to noise in the images and error in the force measurement. Subsequently, we discuss how the speckle size impacts the resolution of inferred features. This leads us to show how the forcing condition also influences the inversion. To further investigate the nuances of $\infty$--IDIC, we introduce an expected information gain (EIG) heuristic that provides comparison between inversions. We use this to further investigate the influence of speckle size and magnitude of applied force. 

\subsubsection{Image Noise and Force Measurement Error}

\begin{figure}[H]
    \centering
    \includegraphics[width=.65\textwidth]{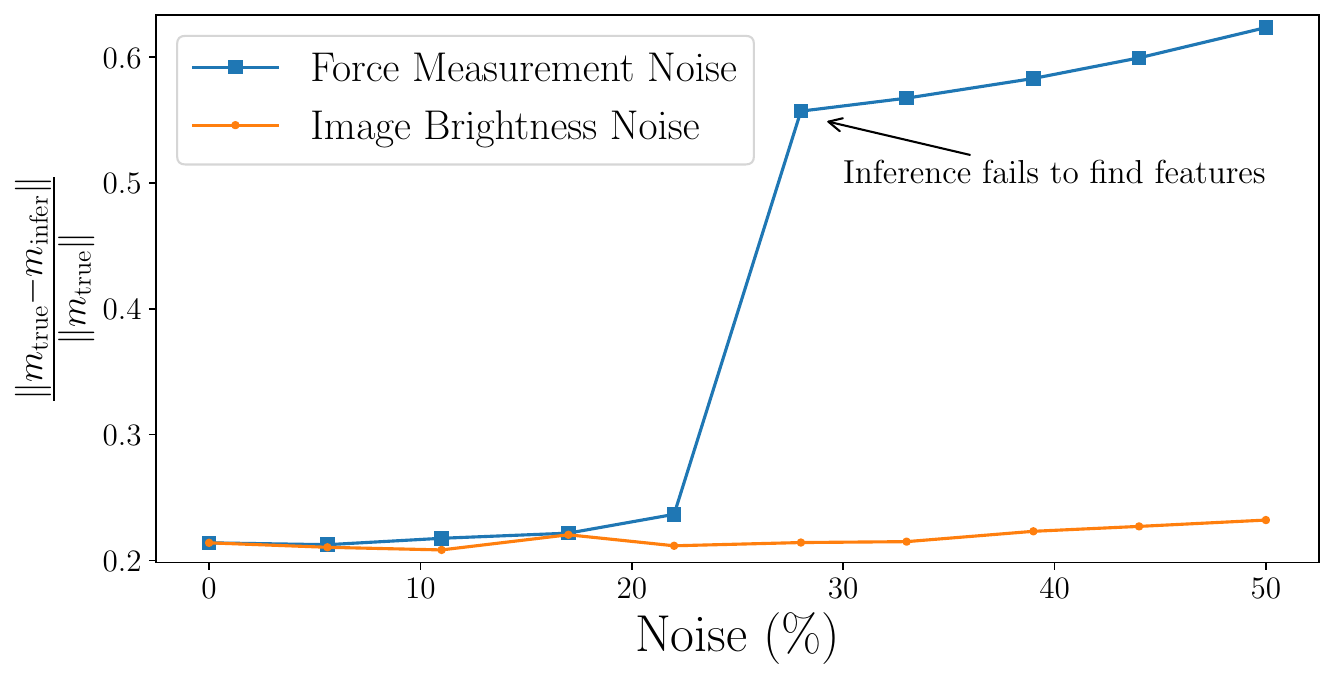}
    \caption{Effect of noise on the accuracy of the inferred results for Hyperelasticity where $m_{\text{true}}$ was a $5 \times 5$ case. Error in the force measurement eventually led to the inferred solutions not discovering the features, whereas image brightness noise had little effect.}
    \label{fig:noise_results_hyperelasticity}
\end{figure}

In real-world experiments observational noise and error are inevitable, so we investigate their influence  on the inferred results. \ref{appendix:a} explains synthetic data generation with noise corruption. We show in Figure \ref{fig:noise_results_hyperelasticity} that with image brightness noise, $\infty$--IDIC infers the correct solution as at 50\% brightness noise and the inversion maintains 20 \% accuracy. Previous studies show integrated DIC inverting for two parameters at 3\% noise in the image brightness values begins to show significant error for single experiments, and requires multiple experiments to achieve the same accuracy \cite{Mathieu2015EstimationIntegrated-DIC}. On the other hand, for error in the force measurement, above 25\% error the inversion fails to find the spatial distribution of the features. A notable advantage of $\infty$-IDIC is its robustness to force measurement errors, which is particularly beneficial in contexts where measuring forces is challenging. For instance, when characterizing composite microstructure inside a Scanning Electron Microscope with IDIC, the force is applied outside the image field of view, necessitating an approximation of the local boundary conditions \cite{Rokos2023} such that being robust to force measurement is helpful. Regardless, we suspect that in most settings, the force measurement error is sufficiently less than 25\%. 

\subsubsection{Relationship Between Speckle Size and Resolution of Inferred Features}

The speckle pattern is used during DIC to provide a contrast for comparison between the reference and deformed images. There are many methods to generate speckle patterns (e.g., spray paint, airbrush, stamping, etc.) and it is of interest to understand how the speckle pattern influences the inversion in high dimensions, particularly when the feature size is on similar length scales to the speckle size. As presented in the appendix, we generate the synthetic speckle by applying a threshold to a Gaussian random field with a tunable correlation length \cite{Villa2021}. We therefore study the effect of speckle size by considering an inverse problem for a parameter field consisting of a single void in the domain, and compare the solutions of the inverse problems for various combinations of speckle and feature sizes.

\begin{figure} [H]
    \centering
    \includegraphics[width=.60\textwidth]{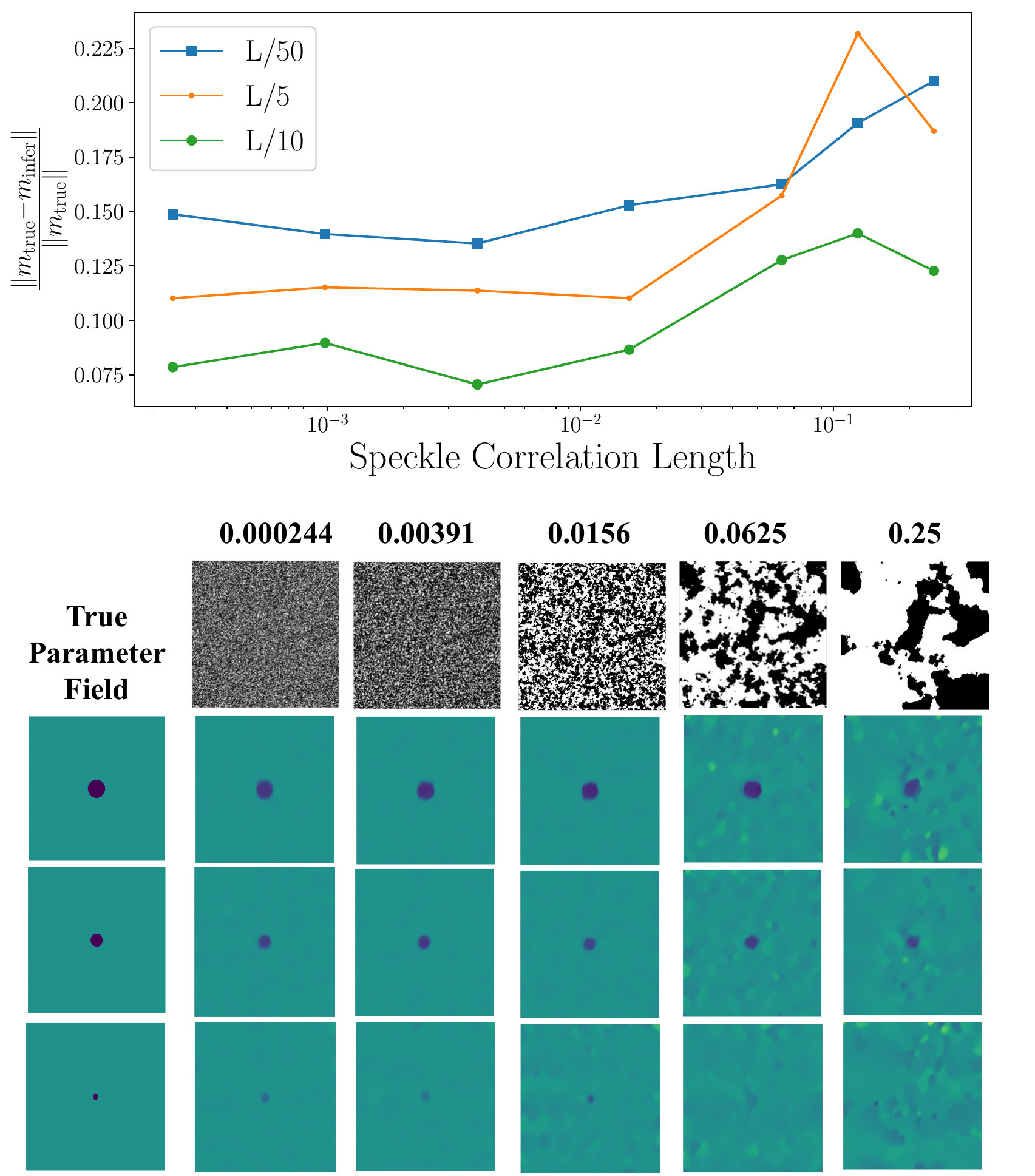}
    \caption{The resolution of inverted features is influenced by the speckle size. $\infty$--IDIC infers the correct solution even when the speckle size is sufficiently smaller than the feature size. 10\% and 5\% noise were applied to the image brightness and force measurements, respectively. The loading condition is a simple tensile experiment with $\sim$5\% strain.}
    \label{fig:speckles_hyperelasticity}
\end{figure}

Figure \ref{fig:speckles_hyperelasticity} shows the inverted parameter fields for the void, where the size of the void (feature) is reduced from $L/5$ to $L/50$ ($L=1$ is the length of the domain) and the speckle correlations are reduced from $2.5 \times 10^{-1}$ to $2.44 \times 10^{-4}$. Generally, we observe that smaller speckle sizes improve inversion results, allowing us to recover smaller features in the underlying parameter field. But, surprisingly, even when the speckle size is larger than the feature it is still discovered. We speculate that the edges of the speckle pattern (where $\nabla I_1$ is largest) contain vital information for $\infty$--IDIC, being the main contribution to the derivative of the misfit functional (see Equation \eqref{eq:misfit_gradient}), such that large speckle sizes can be informative. Yet, we observe that when the feature size is reduced to $L/50$, larger speckle sizes fail to detect the void in the parameter field, while the smallest speckle size can still accurately recover the void. 

\subsubsection{Varying the Forcing Condition} \label{sec:loading_condition}

\begin{figure} [H]
    \centering
    \includegraphics[width=.75\textwidth]{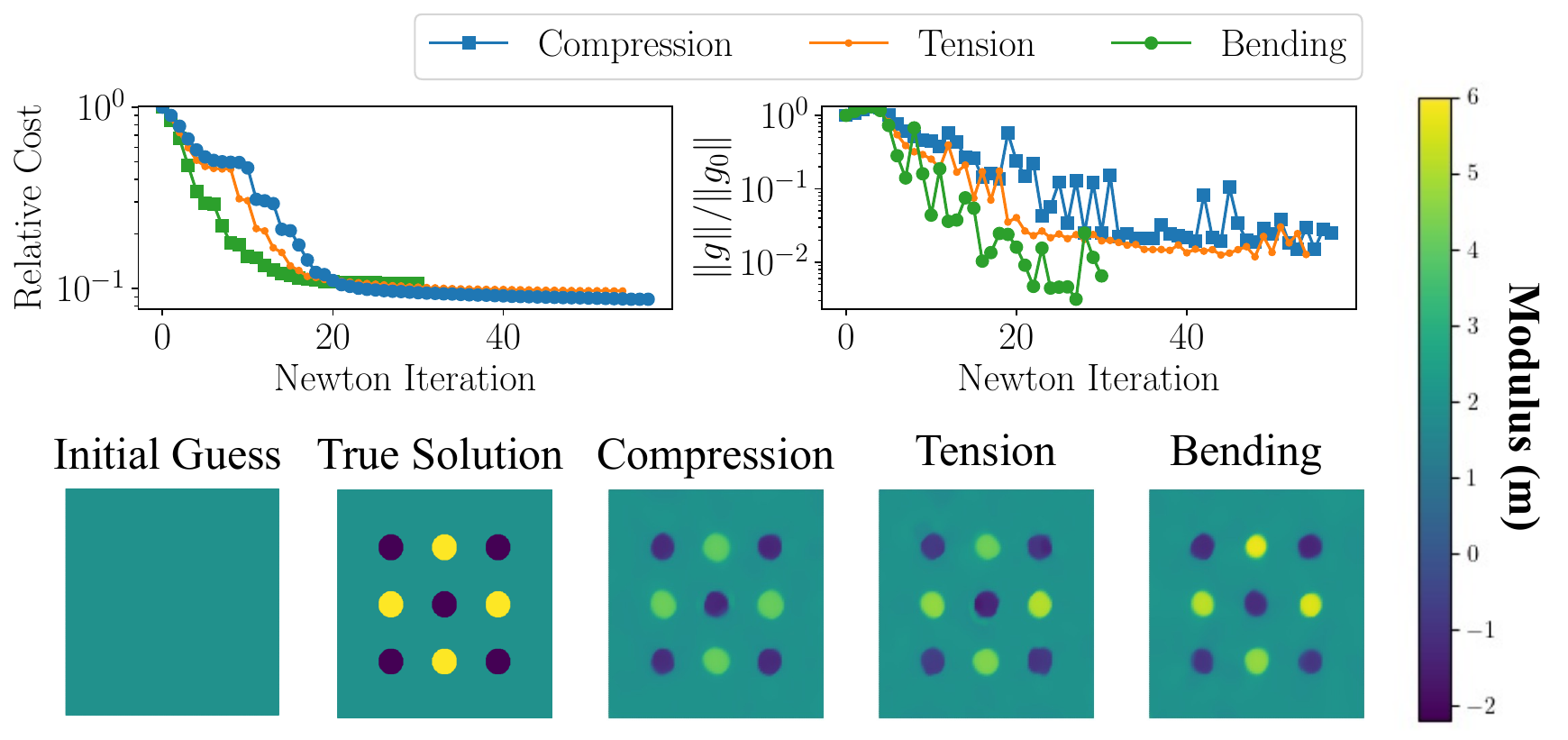}
    \caption{$\infty$--IDIC infers the correct solution for varying experiment boundary conditions. 10\% and 5\% noise were applied to the image brightness and force measurements, respectively. For each of the loading experiments, $\sim$5\% strain was applied. The inverse problem is solved using the primal-dual total variation formulation ($\gamma_{L^2} = 5\times10^{-6}$, $\gamma_{TV} = 7.5 \times10^{-4}$) with a mesh size of $100 \times 100$, and the speckle correlation length is 0.01.}
    \label{fig:boundary_conditions_hyperelasticity}
\end{figure}

In addition to the speckle pattern, the load configuration is a necessary component of the experiment. Figure \ref{fig:boundary_conditions_hyperelasticity} shows inversion solutions of the $3 \times 3$ grid modulus pattern using three different loading conditions; compression ($t_{\text{normal}}$ = -0.50), tension ($t_{\text{normal}}$ = 0.50), and bending ($t_{\text{shear}}$ = -0.25). 
The errors compared to the true solution, $\|m_{\mathrm{true}} - m_{\mathrm{infer}}\|/\|m_{\mathrm{true}}\|$, are 0.233, 0.219, and 0.218 for compression, tension and bending, respectively. We observe that the relative cost and gradient norm reduce more slowly for compression than the other loading conditions. Bending converges the quickest as shown in Figure \ref{fig:boundary_conditions_hyperelasticity} and  appears to provide maximum information at the top of the domain where the sample is in tension. This aligns with the observation that tension is more informative than compression. Section \ref{sec:multi_exp} will consider coupling multiple experiments in a single inverse problem to further improve the inversion by providing more information.

\subsubsection{An Information Gain Heuristic}\label{section:numerics_eig}

We now introduce an expected information gain (EIG) metric from Bayesian optimal experimental design, as a means of evaluating how well-informed different experimental setups are relative to each other. EIG is a statistical quantity that measures how much information is gained about a random variable (in our case $m$) from the observation of another (possibly correlated) random variable (in our case the observation $I_1$). Given a means of computing EIG, we can make relative information comparisons by evaluating EIG for differing experimental setups (e.g., different traction conditions or different speckle patterns). While EIG formally requires a statistical inference framework for its derivation, in the case of linear inverse problems with Gaussian priors, its can be computed from a closed form expression \cite{alexanderian2016bayesian}. The Gaussian prior plays a role analogous to the quadratic regularization (e.g., $H^1(\Omega)$ but not TV) in our presentation, so for this reason we restrict our attention to $H^1(\Omega)$ regularization in this section.
In the context of linear Bayesian inverse problems, EIG is computed as follows. First, one computes the maximum a posteriori (MAP) point, which in this setting coincides with the solution of the deterministic inverse problem $m^*$. Second, one computes the generalized eigenpairs $\{(\lambda_i,\psi_i)|\lambda_i \geq \lambda_j \forall i < j\}_{i=1}^\infty$ at $m^*$, here written in variational form
\begin{equation}
    \left\langle D_{mm}\Phi(u(m^*,t);I_0,I_1) \psi_i , \tilde{m} \right\rangle_\mathcal{M} = \lambda_i \left\langle \mathcal{C}^{-1} \psi_i , \tilde{m} \right\rangle_\mathcal{M}  \qquad \forall \tilde{m} \in \mathcal{M} \label{eq:gevp},
\end{equation}
where $\mathcal{C}^{-1}$ is the self adjoint operator that induces the norm that is used in the regularization. For $H^1(\Omega)$ regularization this operator is the Laplace operator $\mathcal{C}^{-1} = -\gamma \Delta + \delta \mathcal{I}_{\mathcal{M}}$ with homogeneous Neumann boundary conditions, for appropriate choices of $\gamma,\delta$. Then EIG, hereby denoted as $\Psi$ is computed from the following formula
\begin{equation}
    \Psi = \log \det\left(I_\mathcal{M} + \mathcal{C}^{1/2}D_{mm}\Phi \mathcal{C}^{1/2}\right) = \sum_{i=1}^\infty \log(1+\lambda_i) \approx \sum_{i=1}^r \log(1+\lambda_i),\label{eq:eig_formula}.
\end{equation}
where $r$ is chosen such that $\lambda_j \approx 0$ for all $j >r$. 
\begin{remark}
We note that the use of $H^1(\Omega)$ in the context of infinite-dimensional Bayesian inverse problems is abusive, as in this case $\mathcal{C}^{-1}$ is not trace class for $\Omega \subset \mathbb{R}^d$ with $d>1$; thus there is no meaningful limit of this problem. However, for a fixed discretization of the problem, this formulation leads to a meaningful approximation of EIG. For a more detailed discussion of EIG and the Bayesian interpretation of $\infty$--IDIC see \ref{appendix:bayes}.
\end{remark}
In this interpretation of EIG, \eqref{eq:eig_formula} suggests that the curvature of the Hessian located at the solution of the deterministic inverse problems indicates information gained in the solution of the optimization problem. Steeper basins around $m^*$ indicate a well-informed solution of the inverse problem, while flatter basins indicate less informed solutions; this makes sense as in the latter case a perturbation of $m^*$ would lead to a comparatively small change in $\Phi$. Since we are equipped to compute all of these quantities via their derivations in section \ref{subsection:derivative_derivations}, we propose using this linearized form of EIG as an information gain metric to compare trade offs in how the relative information gained between different choices of hyperparameter used in the inverse problem formulation. We note that this EIG metric can be efficiently approximated by the use of randomized generalized eigenvalue solvers \cite{HalkoMartinssonTropp2011,MartinssonTropp2020}, since only the first $r$ eigenpairs are needed as in \eqref{eq:eig_formula}. These eigenvalue solvers are implemented in \texttt{hIPPYlib} \cite{Villa2021}.

\subsubsection{Information Gain Increases with Refining Speckle}

\begin{figure} [H]
    \centering
    \includegraphics[width=.50\textwidth]{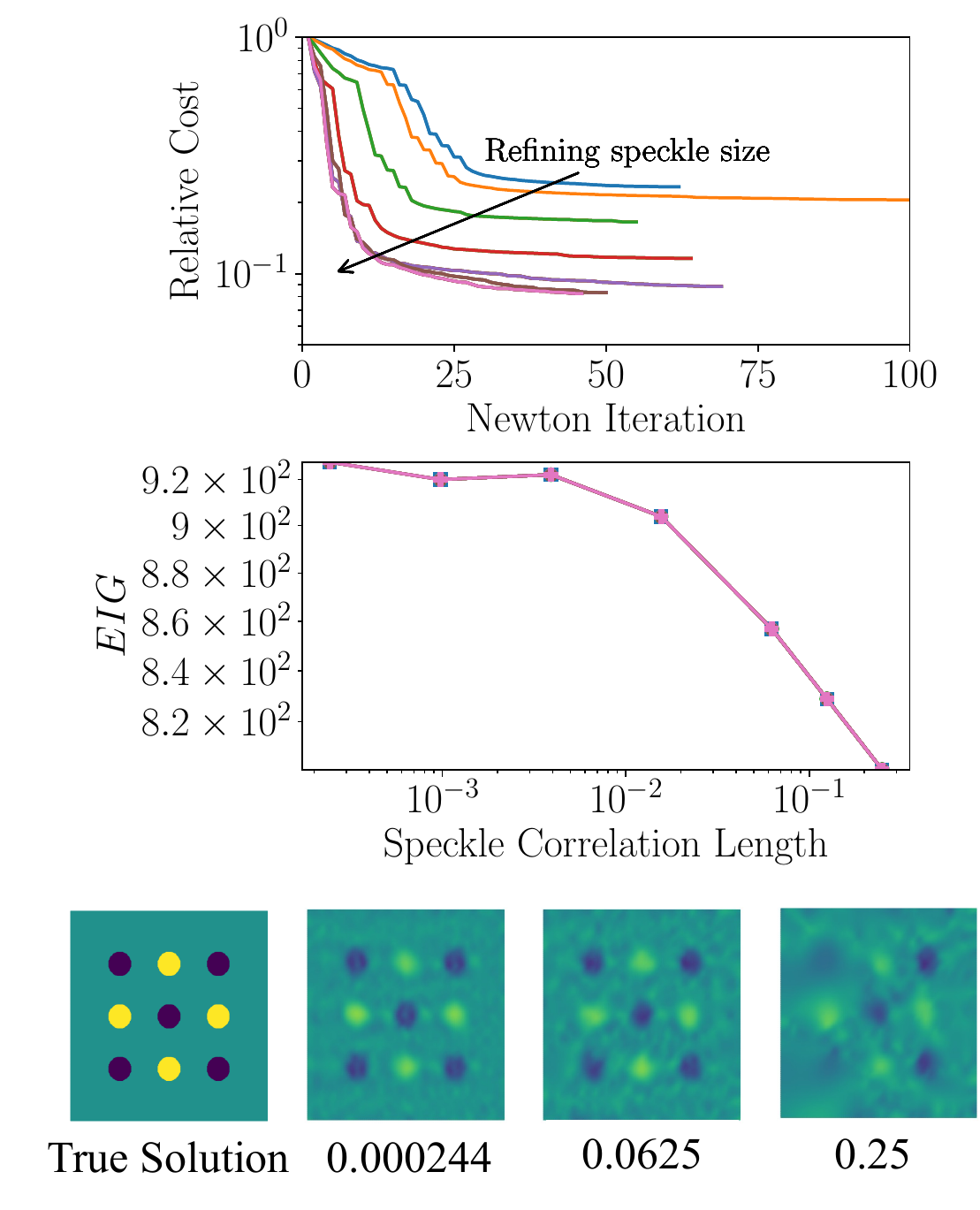}
    \caption{The speckle size of the synthetic images is varied to demonstrate the relationship between the EIG and speckle size. The number of features is held constant at $3 \times 3$ while the speckle size is varied. Here, 0\% noise is applied in order to not corrupt the eigenvalues. The loading condition is a simple tensile experiment with ~5\% strain. An $H^1$ regularization is used.}
    \label{fig:speckles_eig}
\end{figure}

Now that we have defined a heuristic to compare between experiments, we run an experiment where the mesh, number of features, and regularization are held constant at $100 \times 100$, $3 \times 3$ and $\gamma_{H^1} = 5\times10^{-4}$, while the speckle correlation length is refined. Figure \ref{fig:speckles_eig} demonstrates that the EIG improves until plateauing around a speckle correlation length of 0.0156. This aligns with the observation in Figure \ref{fig:speckles_hyperelasticity} that the inversion improves with refining speckle size, until the speckle is sufficient for the feature size. Notably, the relative cost in Figure \ref{fig:speckles_eig} also decays quicker with increasing speckle size suggesting that less Newton iterations are required. The relative cost, EIG, and inverted parameter fields are in agreement that refining the speckle size leads to improved inversions until a certain point where there is a diminishing return. 

\subsubsection{Information Gain Increases with Larger Forcing} \label{sec:loading_EIG}

\begin{figure} [H]
    \centering
    \includegraphics[width=.5\textwidth]{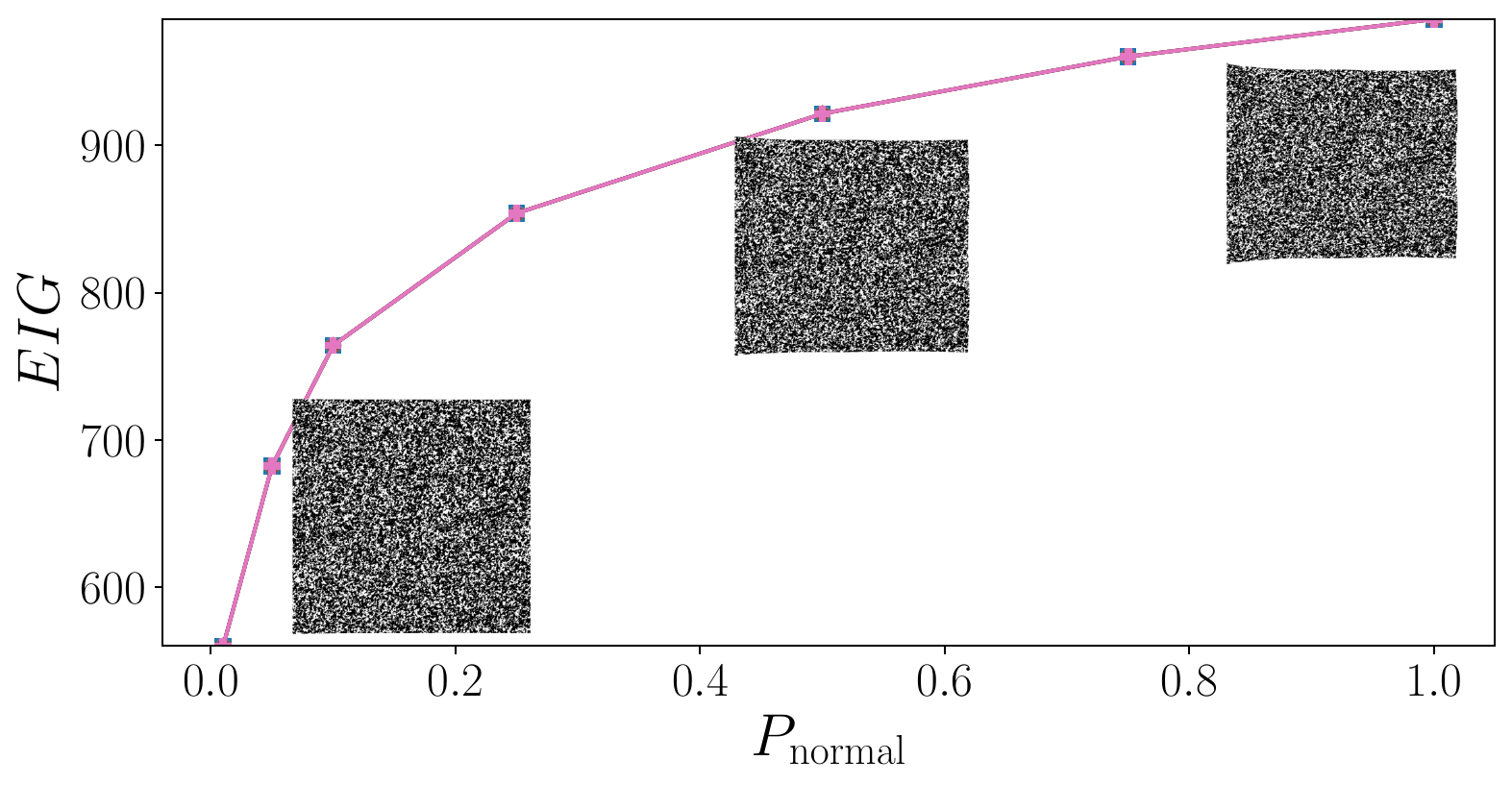}
    \caption{The applied normal force, $t_{\text{normal}}$ of the synthetic images is increased to demonstrate that additional deformation improves the inversion. The number of features is held constant at $3 \times 3$ while the speckle size is varied. Here, 0\% noise is applied in order to not corrupt the eigenvalues. The loading condition is a simple tensile experiment with ~5\% strain. An $H^1$ regularization is used.}
    \label{fig:force_EIG}
\end{figure}

Recall that an important aspect of DIC experiments to consider is the applied forces, as the loading condition determines the displacement observations. We again consider the $3 \times 3$ grid modulus pattern, but now we solve the inverse problems corresponding to various tensile loads from $t_{\mathrm{normal}} = 0$ to $t_{\mathrm{normal}} = 1$. As in the EIG analysis for speckle size, 
we use 0\% noise for the synthetic images and adopt an $H^1$ regularization in the inversion. Notably, often in DIC, large deformations are difficult to handle due to ill-posedness and so a series of images to take smaller steps \cite{grant2009high}. Yet, since IDIC benefits from incorporating the governing equations of the continuum mechanics we observe that a single pair of images can be used at large deformations ($\sim$ 10\%) as seen in Figure \ref{fig:force_EIG}. We plot the resulting EIG computed at the MAP points also in Figure \ref{fig:force_EIG}. We observe that as the applied force increases, the EIG increases and $\infty$--IDIC is able to resolve the modulus field better. By increasing the deformation of the speckle pattern, $\infty$--IDIC is provided with more informative data. 


\subsection{Useful Extensions of $\infty$--IDIC}

We will present two forward-looking aspects of $\infty$--IDIC. First, we consider incorporating multiple experiments into a single inverse problem. Then, we show that we can discover stress concentrations by capturing the effects of the spatially-varying material properties.

\subsubsection{Learning More Through Multiple Experiments} \label{sec:multi_exp}

We showed that increasing the magnitude of the force improved the EIG (see section \ref{sec:loading_EIG}) and the inversion solution changed varying the loading configuration (see section \ref{sec:loading_condition}). This section demonstrates the marriage of multiple experiments into a single inverse problem; It has been shown before that using multiple experiments improved the inversion of five isotropic elastic parameters in an IDIC formulation \cite{Neggers2019SimultaneousIdentification}. To formulate $\infty$--IDIC with multiple experiments, we update the cost functional for $\infty$--IDIC to be a sum over the experiments,

\begin{equation} \label{eq:cost_multiple_experiments}
	\widehat{\mathcal{J}} (m) := \left( \frac{1}{N_{\text{exp}}} \sum_{i=1}^{N_{\text{exp}}} \Phi \left(u(m, t_i); I_{0,i},I_{1,i} \right) \right) + \mathcal{R}(m),
\end{equation}
where $N_{\text{exp}}$ is the number of experiments, $(t_i)_{i=1}^{N_{\text{exp}}}$ are the different loading conditions, and $(I_{0,i})_{i=1}^{N_{\text{exp}}}$ and $(I_{1,i})_{i=}^{N_{\mathrm{exp}}}$ are the undeformed and deformed images corresponding to each experiment.
Similarly, the forward, adjoint, and gradient equations become a sum over the experiments and subsequently the incremental forward, incremental adjoint, and incremental gradient as well. Including the $1/N_{\text{exp}}$ ensures that the regularization weighting is consistent with a single experiment. 

We first show results for a single tension ($t_{\text{normal}}$ = 0.50) experiment to contrast handling four experiments simultaneously: tension ($t_{\text{normal}}$ = 0.50), compression ($t_{\text{normal}}$ = -0.50), downwards bending ($t_{\text{shear}}$ = -0.10), and upwards bending ($t_{\text{shear}}$ = 0.10). The $3 \times 3$ feature example noticeably better recovers the true solution's modulus values with the additional experiments. Small features that were not recovered with a single experiment are improved such as the ears of the Texas Longhorn. Lastly, spurious features are removed upon with the addition of more experiments. This is especially the case with the Voronoi tessellations. We find that varying the loading conditions improves the inversion for $\infty$--IDIC. 

\begin{figure} [H]
    \centering
    \includegraphics[width=.60\textwidth]{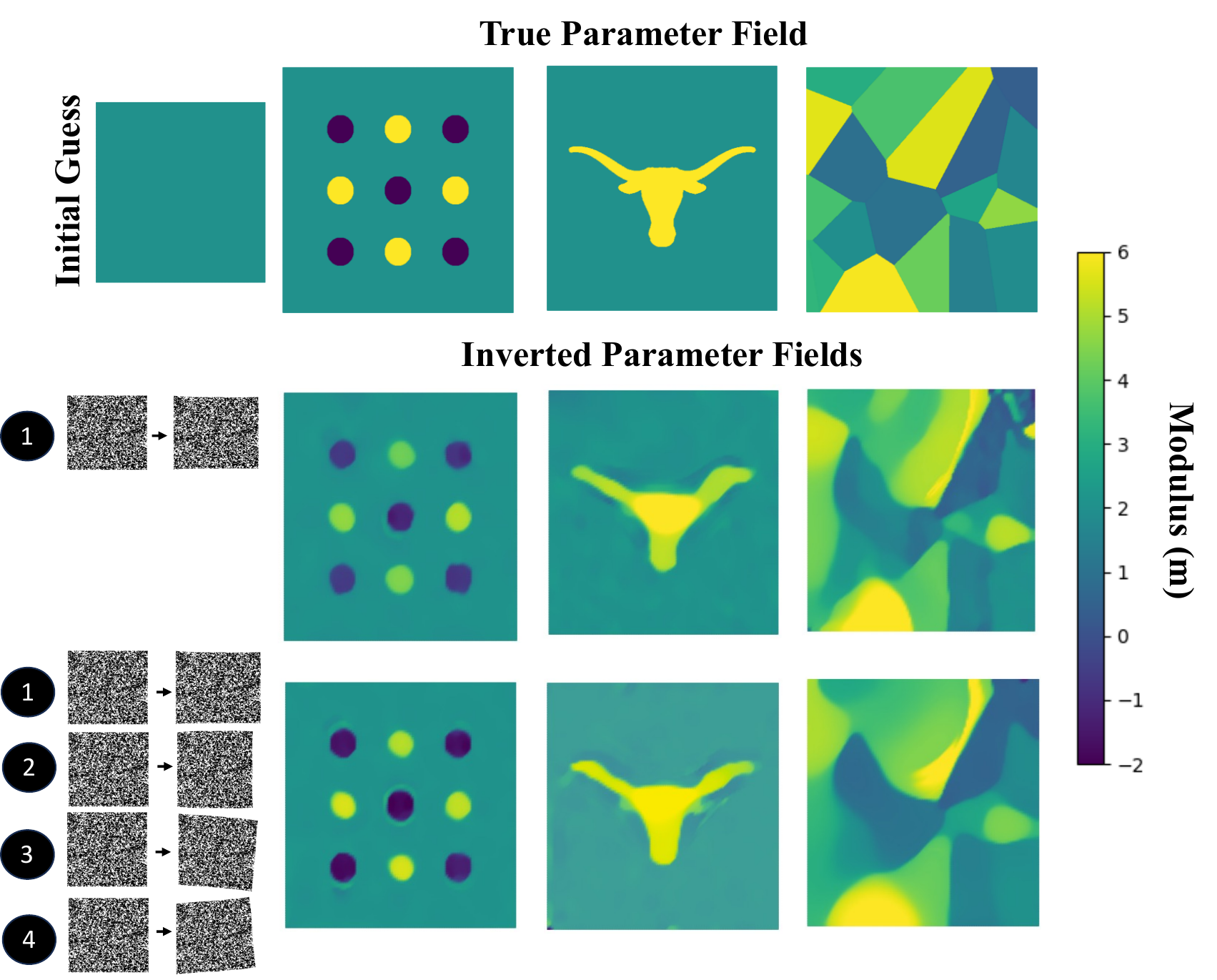}
    \caption{$\infty$--IDIC infers the correct solution for multiple experiments. 10\% and 5\% noise were applied to the image brightness and force measurements, respectively. For each of the loading experiments, $\sim$5\% strain was applied, albeit $t_{\text{normal}}$ and $t_{\text{normal}}$ varied. The inverse problem is solved using the primal-dual total variation formulation ($\gamma_{L^2}$ = 5e-6, $\gamma_{TV}$ = 7.5e-4) with a mesh size of $100 \times 100$, and the speckle correlation length is 0.01.}
    \label{fig:multiple_experiments_hyperelasticity}
\end{figure}

\subsubsection{Post-processing for Modulus-Induced Stress Concentrations}

In DIC experiments, it is customary to post-process the displacement field to a strain field using the Green--Lagrange strain relationship, 

\begin{equation} \label{eq:strain}
    \mathbf{\epsilon} = \frac{1}{2} \left( \nabla u + \nabla u^T \right).
\end{equation}

\begin{figure}[H]
    \begin{minipage}[c]{0.5\textwidth}
        \centering
        \includegraphics[width=\textwidth]{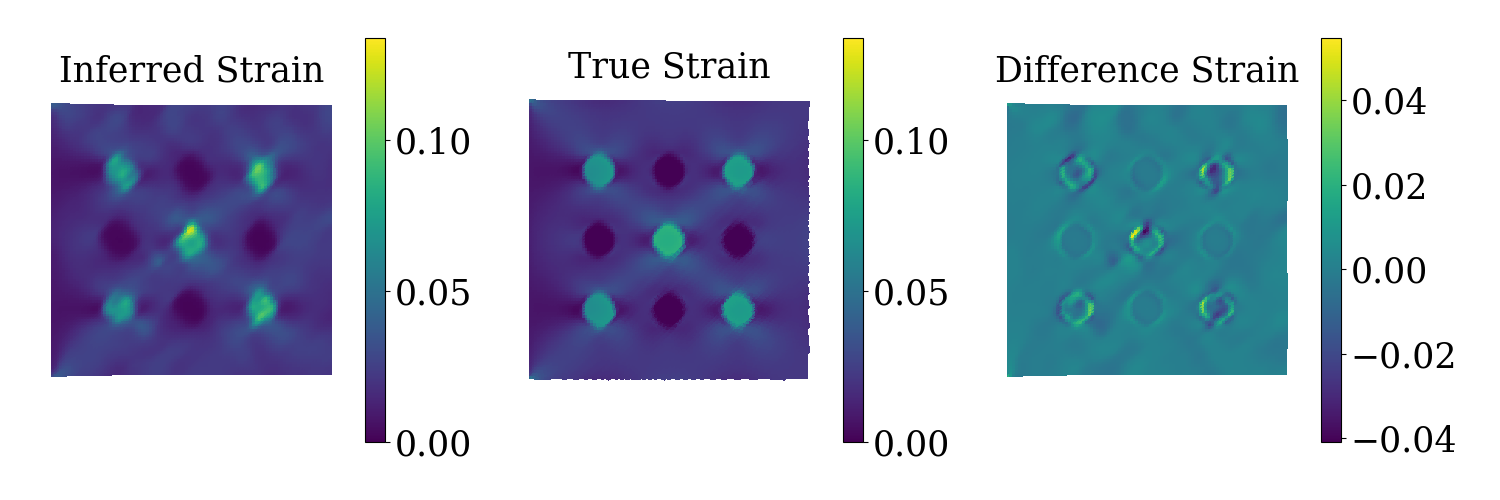}
        \includegraphics[width=\textwidth]{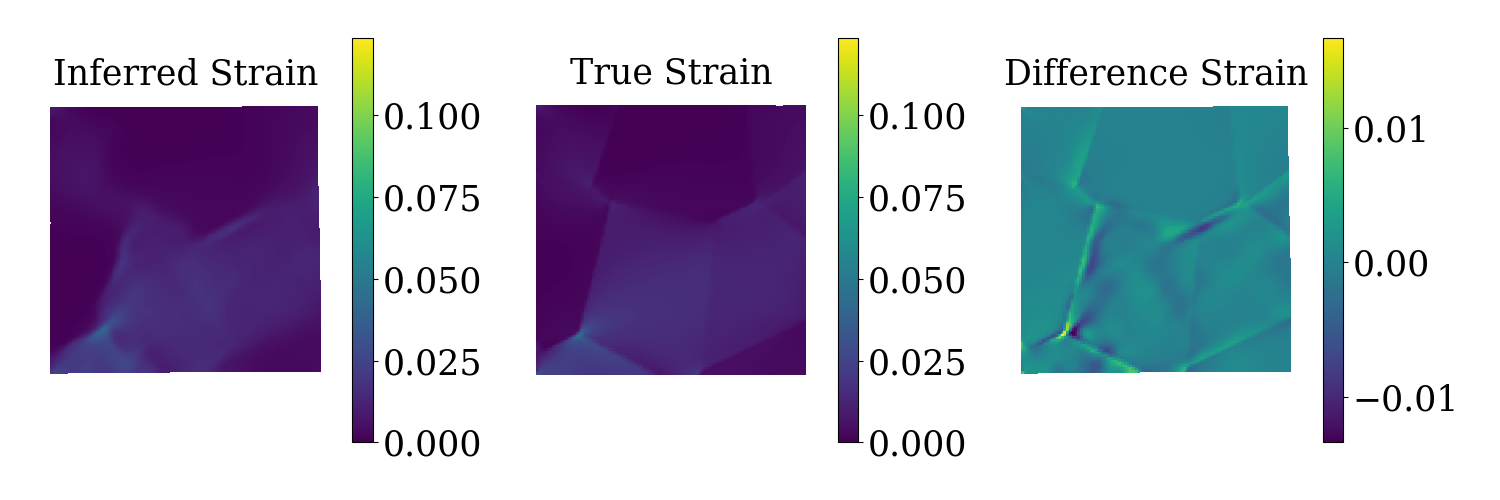}
        \includegraphics[width=\textwidth]{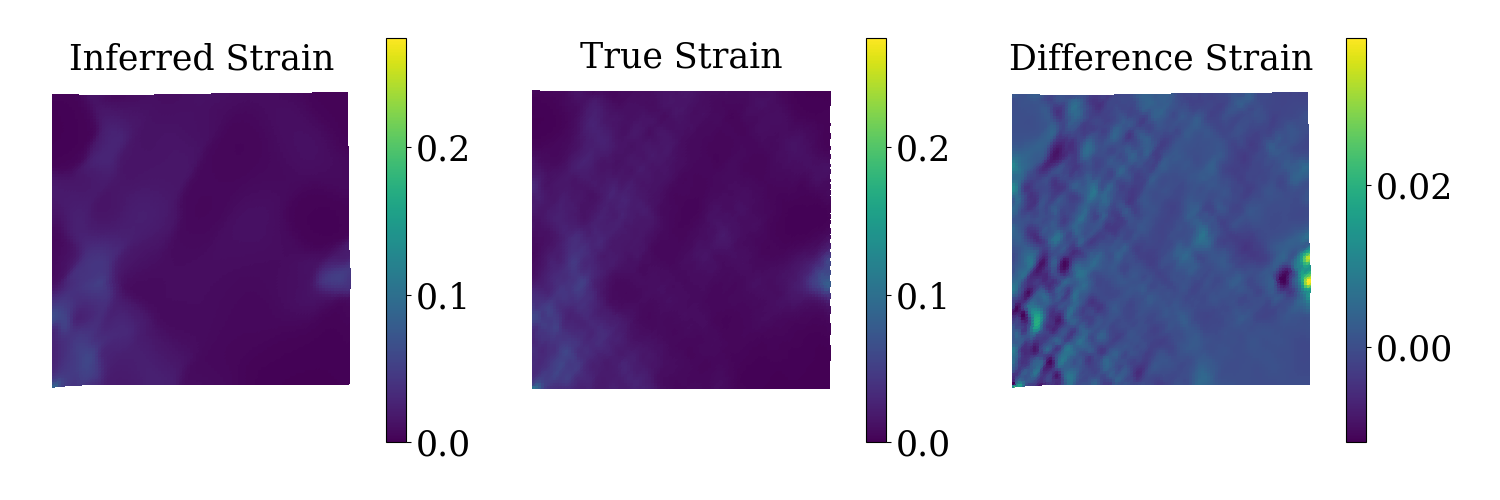}
    \end{minipage}
    \begin{minipage}[c]{0.5\textwidth}
        \centering
        \includegraphics[width=\textwidth]{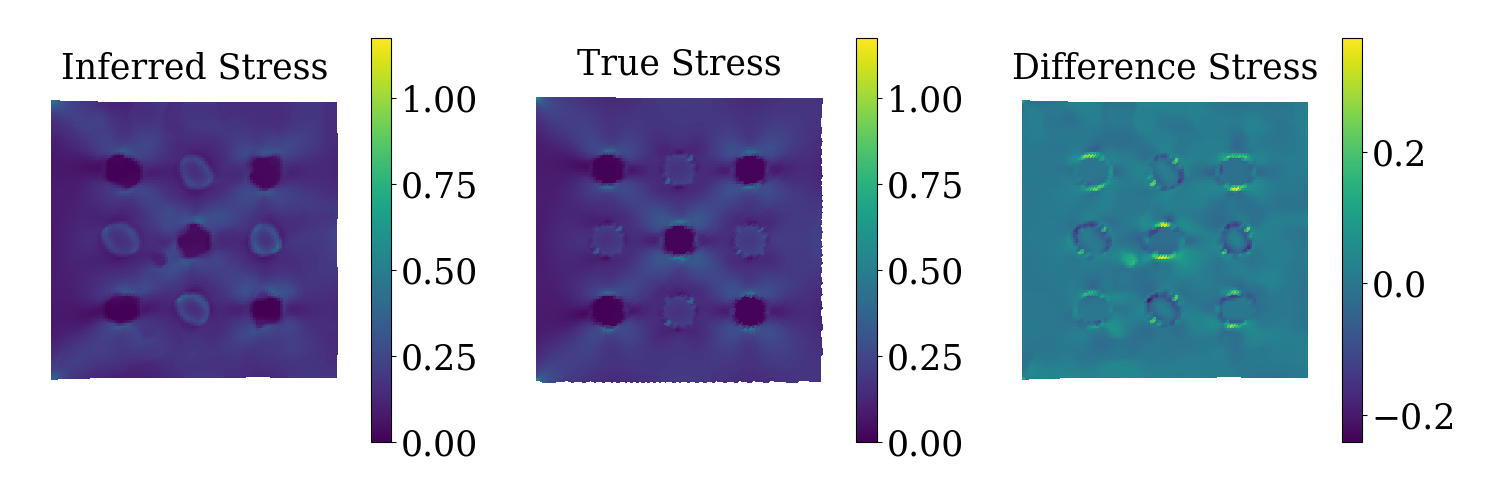}
        \includegraphics[width=\textwidth]{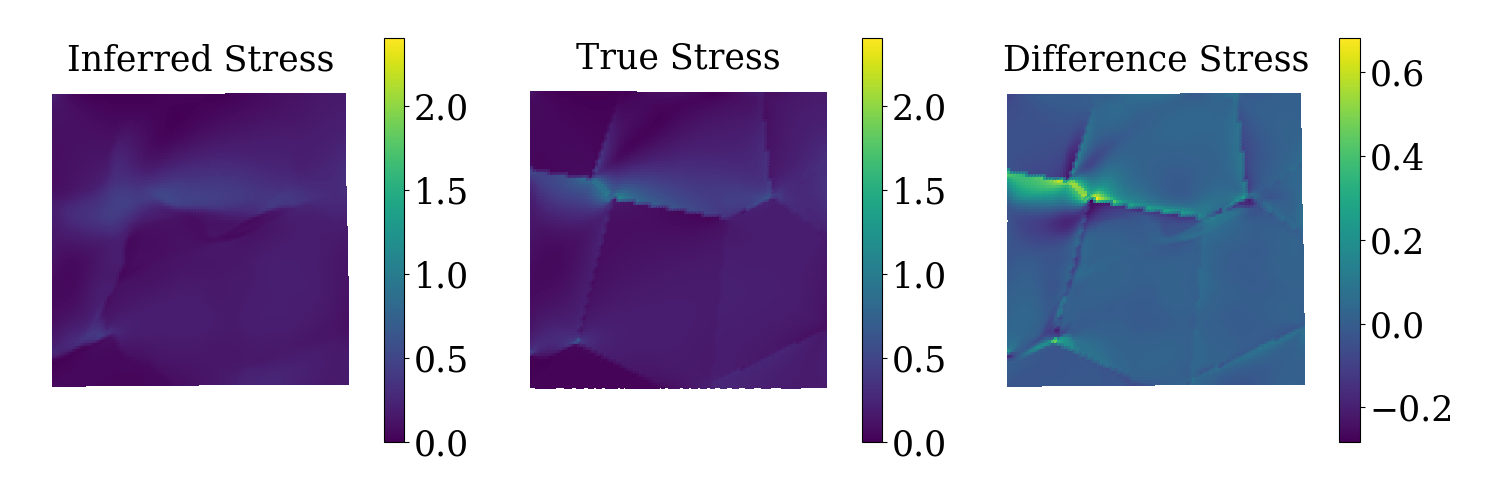}
        \includegraphics[width=\textwidth]{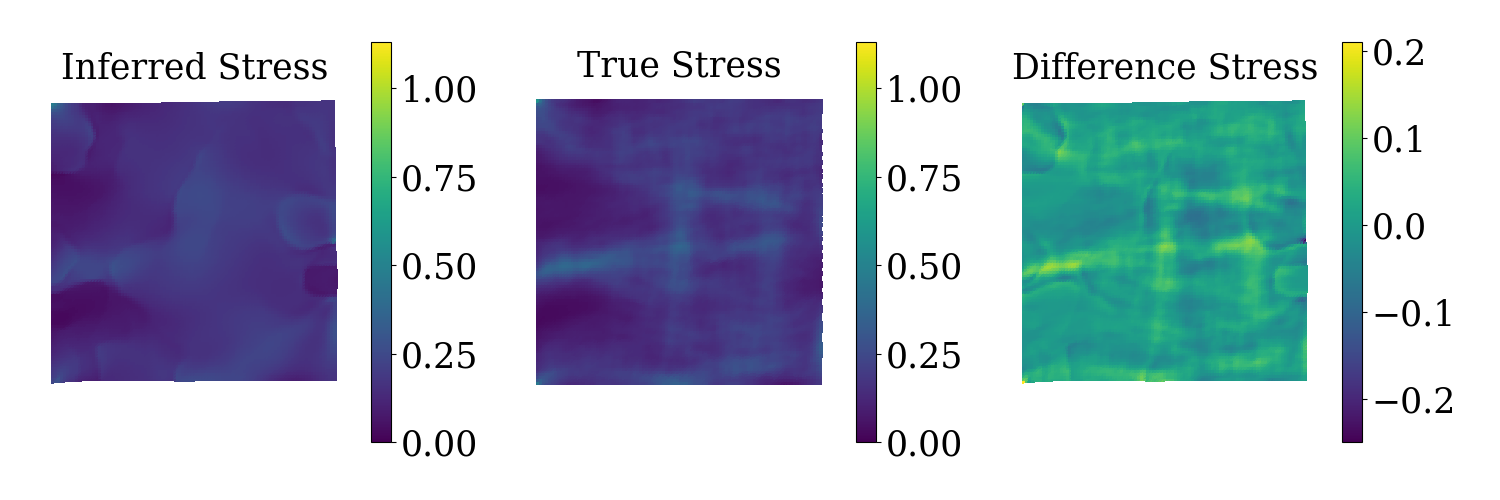}
    \end{minipage}
    \caption{Strain and stress fields for the true and inferred solutions for linear elasticity. The first row shows the strain fields for the $3 \times 3$ grid of alternating voids and stiff regions. The second row shows the strain fields for the Voronoi tessellation. The third row shows the strain fields for the isotropic Gaussian random field.}
    \label{fig:inverse_results_linear_2}
\end{figure}

For instance, Ncorr, an open-source DIC software, provides this functionality \cite{Blaber}. The strains are visualized as spatially varying scalar fields by selecting one direction in the tensor, such as $\varepsilon_{xx}$ which aligns with the directional of the normal applied load. We show results of $\varepsilon_{xx}$ using linear elasticity shown in Figure \ref{fig:inverse_results_linear_2}. The strain fields provide insight into loading behavior. Yet, engineering systems are often designed to stress criteria. In the setting of heterogeneous materials, the stress field is \emph{different} from the strain field. We use a common stress criterion for linear elasticity, the von Mises stress which is used to predict yielding in materials under uniaxial tension. The von Mises stress is a scalar field that represents the stress state of the material and is convenient for visualization. For linear elasticity, a von Mises stress field, $\sigma_{vm}$, can be calculated as

\begin{subequations}
\begin{alignat}{1}
\label{eq:vm_stress}
        &\sigma_{vm} = \sqrt{\frac{3}{2} s : s} \\
\label{eq:s_stress}
        &s = \sigma - \frac{1}{3} \text{tr}(\sigma) I. 
\end{alignat}
\end{subequations}

We show post-processing results for linear elasticity only, but similar techniques could be done for the hyperelasticity case. Figure \ref{fig:inverse_results_linear_2} shows the true and inverted strain and stress fields for the same features as in Figure \ref{fig:inverse_results_linear}. The first example involves a $3 \times 3$ grid of alternating voids stiff regions where the inferred strain and stress fields capture the general trends successfully. The strain field shows limited deformation in stiff interfaces and large deformation in the void regions while stress concentrations are observed around voids, indicating potential high-risk regions for failure. In the Voronoi tessellation, the strain field shows maximum strain in grains with lower modulus and minimal strain in high modulus regions. The inferred fields properly identify stress concentrates at grain interfaces, aligning with the true solution. Lastly, with the isotropic Gaussian random field there are stress concentrations properly identified around high modulus regions. $\infty$--IDIC captures general trends in strain and stress fields even for a single experiment at low strains ($\sim$0.2\%), leading to potential applications in non-destructive failure predictions. 

\section{Conclusion}

In summary, we present $\infty$--IDIC, a scalable and efficient framework for the inference of spatially-varying heterogeneous materials that parametrize continuum mechanical models. The goal of $\infty$--IDIC is the characterization of damage, defects, anomalies and inclusions in material specimen through nondestructive evaluation methods, without any prior knowledge of their morphology or strength. These effects cannot be directly observed, but must be reconstructed indirectly through inverse problems. Since heterogeneous material fields are formally infinite dimensional, we require (i) rich information in the form of observational data to faithfully reconstruct the material properties, and (ii) an infinite-dimensional formulation of the inverse problem. To address the first concern we build on IDIC, which like traditional DIC utilizes images of material deformation as observational data, but unlike DIC monolithically couples the image registration to the governing equations of equilibrium to enforce physical consistency. To address the second concern we formulate a general formulation in a function space setting to allow for the inference of arbitrary material parameter fields. The main challenges of this approach are due to the mathematical and computational challenges posed by the infinite-dimensional formulation, namely ill-posedness and high dimensionality; addressing these challenges is a focus of this work. First, to tackle the inherent ill-posedness we consider various regularization schemes, namely $H^1(\Omega)$ and total variation for the inference of smooth and sharp features, respectively. To address the computational costs associated with the discretized problem, we utilize an cost efficient and dimension-independent inexact-Newton CG framework for the solution of the regularized inverse problem.

In the numerical results, we demonstrate that $\infty$--IDIC can recover spatially varying material properties with remarkable accuracy and sharpness using examples of linear elasticity and hyperelasticity. We demonstrate the ability to invert for 40,401, 160,801, and 361,201 parameter values while maintaining mesh-independent behavior. We conduct numerical experiments to investigate how robust $\infty$--IDIC is to noise in the image and force data, variation in speckle sizes, and varying force conditions. We find that refining the speckle size and increasing the force improves inversions. This leads us to show that coupling multiple experiments in a single $\infty$--IDIC inverse problem further improves inversions. Lastly, we demonstrate that along with the inversion of the spatially varying features, we can invert for stress concentrations arising from heterogeneity which may potentially be useful for non-destructive evaluation.

We recognize that there is a need to characterize complex, heterogeneous fields such as welds, composites, tissues, etc. We believe that $\infty$--IDIC can be a powerful tool for such challenges due to its ability to invert for sharp features and its mesh-independent behavior. We hope that this enables researchers to better characterize complex materials and ultimately build more accurate models. The $\infty$--IDIC framework can be extended to other settings such as fluid mechanics where an image registration problem can be coupled with a PDE model (e.g., using particle image velocimetry observational data). In this work we have focused on 2D PDE problems with simple geometries for the sake of a methodological exposition. Extending $\infty$--IDIC to complex geometries, three-dimensional and time-dependent problems which would be a significant step forward for material characterization, and constitutes future work.

\section{Acknowledgments}
The authors are grateful to Umberto Villa for insightful discussions, and help with software questions related to \texttt{hIPPYlib} \cite{Villa2021}. The authors are also grateful for Graham Pash for insightful discussions.

\newpage

\bibliographystyle{model1-num-names}
\bibliography{Joseph.bib}

\appendix

\section{Details on Numerical Implementation}
\label{appendix:a}

\par
\begin{figure}[H]
\centering
\includegraphics[width=0.5\textwidth]{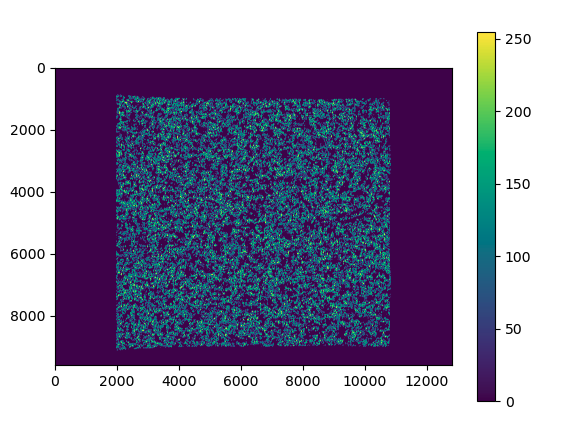}
\caption{Image brightness noise at 50\% added to the synthetically generated images. The noise is added to the images to simulate the effect of varying lighting conditions or general imaging errors in the experiment.}
\end{figure}

This section describes the process of generating synthetic data for the inverse problem. The assumption of a prior agreement between our understanding of the physics and the observed data leads us to commit what is known as an "inverse crime." Each numerical result is derived using the same PDE that is subsequently utilized to inform the inverse problem. Synthetic data is generated by solving the PDEs for the forward problem and saving speckled images before and after deformation. The speckled images are then used as the observational data in the subsequent inverse problem. In an actual experiment, the images would be captured using a camera, or similar sensor, while the speckled specimen is loaded. The following steps detail how we synthetically form speckled image data:

\begin{enumerate}
    \item \textbf{Mesh construction}: Construct a mesh of size $500 \times 500$ where $500$ is the number of elements in the mesh. The mesh is used to approximate $\Omega_0$ via linear Lagrangian triangular elements to form finite element spaces $\mathcal{M}^h \in \mathcal{M}$ and $\mathcal{V}^h \in \mathcal{V}$. The parameter and state space are 1,002,001 and 502,002 dimensions, respectively.
    \item \textbf{True parameter generation}: Generate a true field, $m_{\text{true}}(x)$, that represents the true material model. This field is used to generate the Young's modulus, $E(x) = e^{m_{\text{true}}(x)}$, and the Poisson's ratio $\nu(x) = 0.35$. The true field is generated through a handful of strategies geared towards representing complex material fields. 
    \item \textbf{Boundary conditions definition}: Define the traction condition, $t$, on the right boundary, $\Gamma_{R}$, by determining the forces, $P$, and the Dirichlet boundary condition, $u = 0$, on the left boundary, $\Gamma_{L}$. The Neumann boundary condition is defined as $t = 0$ on the top and bottom boundaries, $\Gamma_{T}$ and $\Gamma_{B}$. Figure \ref{fig:deformation_schematic} illustrates the boundary conditions.
    \item \textbf{Forward PDE solve}: Solve the forward PDE for the true material model to obtain the true displacement field, $u_{\text{true}}(x)$, and the true stress field, $\sigma_{\text{true}}(x)$.
    \item \textbf{Image generation}: Generate a speckled image existing in the discretized parameter space, $\mathcal{M}$ by taking a sample, $i$, from a bilaplacian prior. The correlation length of the bilaplacian prior, $\ell$, regulates the respective correlation length of the speckles. As the speckles in DIC are black (0) or white (255), thresholding is applied. This is done by projecting $i$ to a function space, $\mathcal{M}^h$, and transforming for the image, $I(x)$, using the following equation:
    \begin{equation} \label{eq:speckle_thresholding}
        I(x) = 255 - 255 \, \text{tanh}(100i(x) + 1).
    \end{equation}   
    \item \textbf{Image mapping in deformed states \& Brightness Error}: Map the simulated speckle to both the reference and deformed states, thereby generating the reference image ($I_{0}$) and the deformed image ($I_{1}$). A modified version of the FEniCS plotting tool operating in displacement mode is used. Instead of representing displacement, the color map represents the grayscale speckle values. The displaced image is corrupted with white noise, i.e., $N \sim \mathcal{N}(0, \sigma^2)$, where sigma is the noise level multiplied by the maximum pixel value to simulate image brightness error. Notably, in practice, the material domain, $\Omega_0$, is smaller than the image domain, $\Omega_I$. To simulate this, the images are saved as .png files with additional white space around the material domain.
    \item \textbf{Force measurement error simulation}: To mimic errors in force measurement data, corrupt the traction condition. For instance, if the true force value is $ t_{\text{normal}} = 10$, the corrupted value is $ t_{\text{normal}} = 9.5$ given a $5\%$ error in a tensile measurement. The corrupted value is saved to a .txt file for subsequent use in the inverse problem.
\end{enumerate}

\begin{figure} [H]
    \centering
    \includegraphics[width=\textwidth]{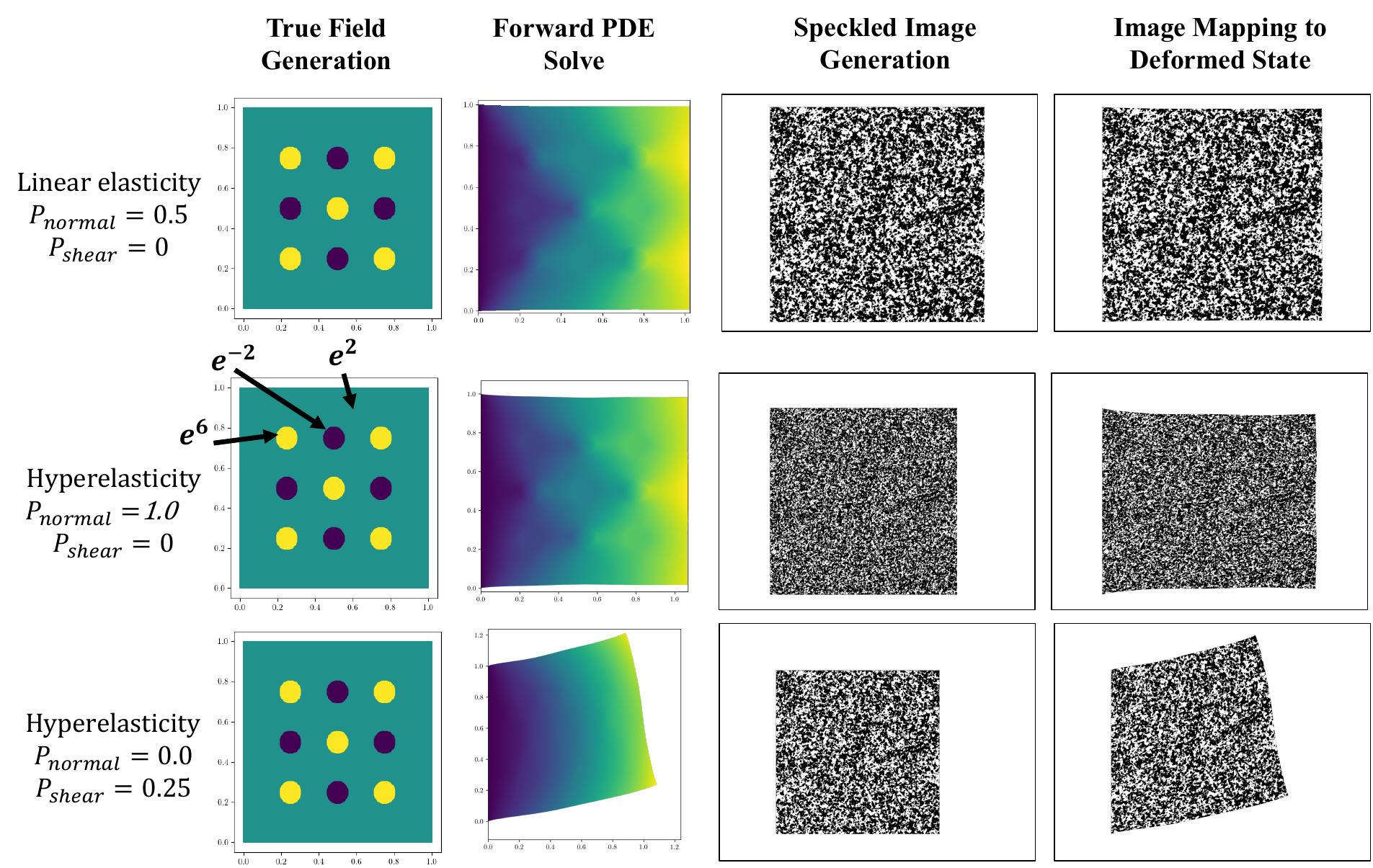}
    \caption{Schematic of the synthetic data generation process.}
    \label{fig:synthetic_generation}
\end{figure}

\section{Bayesian $\infty-$IDIC and expected information gain}\label{appendix:bayes}

We here present a Bayesian extension of $\infty-$IDIC, in order to give sufficient background for the expected information gain (EIG) numerical results presented in \ref{section:numerics_eig}. First we assume $\mathcal{M}$ to be a separable Hilbert space, so basis representations are countable. We denote by $(\mathcal{M},\mathcal{B}(\mathcal{M}))$, the measurable space with $\mathcal{B}(\mathcal{M})$ being the Borel $\sigma-$algebra generated by open sets in $\mathcal{M}$. Next we let $\mathcal{C}:\mathcal{M}\rightarrow \mathcal{M}$ be a trace class, self adjoint linear operator on $\mathcal{M}$, such that
\begin{equation}
    \text{tr}(\mathcal{C}) = \sum_{i=1}^\infty \langle \psi_i,\mathcal{C}\psi_i \rangle_\mathcal{M} < \infty, \qquad \mathcal{C} = \mathcal{C}^*,
\end{equation}
where $\{\psi_i\}_{i=1}^\infty$ is any orthonormal basis for $\mathcal{M}$.
The operator $\mathcal{C}$ can induce an inner product and thereby a norm. Regularization schemes that are quadratics of the norm $\|\cdot\|_{\mathcal{C}^{-1}}^2$ induce Gaussian measures in infinite dimensional settings. The associated regularization 
\begin{equation}
    \mathcal{R}_{\mathcal{C}^{-1}}(m) = \|m - m_0\|^2_{\mathcal{C}^{-1}} = \| \mathcal{C}^{-1/2}(m - \bar{m})\|_{L^2(\Omega)}^2,
\end{equation}
play the same role as a prior distribution, $\mu_\text{prior} = \mathcal{N}(m_0,\mathcal{C})$, in the Bayesian setting. Here $\mathcal{C}^{-1/2}$ denotes the symmetric square root of the operator $\mathcal{C}^{-1}$.

What is left to work towards Bayes is a likelihood measure; towards this end, we assume an additive Gaussian noise model for the observed images, e.g., where noise $\xi \sim \mathcal{N}(0,\mathcal{C}_\text{noise})$, where $\mathcal{C}_\text{noise}:\mathcal{I} \rightarrow \mathcal{I}$ is a self adjoint trace class linear operator. In this case our image misfit takes the following form
\begin{equation}
	\Phi(u(m);I_0,I_1) = \left\|\mathcal{C}^{-\frac{1}{2}}_\text{noise}(I_1(x+u(x)) - I_0(x))\right\|^2_{L^2(\Omega)}
\end{equation}
and leads to the following likelihood
\begin{equation}
	\pi_\text{like}(I_1|m) \propto \text{exp}{\left(-\Phi(u(m);I_0,I_1)\right)}.
\end{equation}

Bayes Theorem provides a rigorous framework for statistically inferring $m$ as a \emph{posterior distribution}, $\mu_\text{post}^{I_1}=\mu_\text{post}(m|I_1)$ given the observational data $I_1$ through the lens of the PDE model:
\begin{equation}
	\frac{d\mu_\text{post}^{I_1} }{d\mu_\text{prior}} \propto \pi_\text{like}(I_1|m).
\end{equation}

Here $\frac{d\mu_\text{post}^{I_1}}{d\mu_\text{prior}}$ denotes the Radon--Nikodym derivative of $\mu_\text{post}^{I_1}$ with respect to $\mu_\text{prior}$. Assuming that the map $m\mapsto u(m)$ is $\mu-$a.e. well-defined, locally Lipschitz continuous and sufficiently bounded the Bayesian inverse problem is well-posed \cite[Corollary 4.4]{Stuart2010}. The expected information gain (EIG) is then defined as the expectation of the Kullback-Leibler divergence ($\mathcal{D}_{KL}$) over instances of data $I_1$ given the noise distribution:
\begin{align}
	\Psi &= \mathbb{E}_{I_1}\left[\mathcal{D}_{KL}(\mu_\text{post}(\cdot|I_1)||\mu_\text{prior}) \right]\\
	\mathcal{D}_{KL}(\mu_\text{post}(m|I_1)||\mu_\text{prior}(m)) &= \int_\mathcal{M} \ln \left( \frac{d\mu_\text{post}^{I_1}}{d\mu_\text{pr}}\right) \mu_\text{post}(dm).
\end{align}
As was discussed in section \ref{section:numerics_eig}, in the case that the forward map $m\mapsto u(m)$ is linear, EIG admits a closed form expression related to the following generalized eigenvalue problem, which is computed at the maximum a posteriori (MAP) point $m^*$,
\begin{subequations}
\begin{align}
m^* &= \text{argmax}_{m \in \mathcal{M}} \mu_\text{post}(m|I_1)\\
\left\langle D_{mm}\Phi(u(m^\star,t);I_0,I_1) \psi_i , \tilde{m} \right\rangle_\mathcal{M} &= \lambda_i \left\langle \mathcal{C}^{-1} \psi_i , \tilde{m} \right\rangle_\mathcal{M}  \qquad \forall \tilde{m} \in \mathcal{M} \\
\Psi &= \sum_{i=1}^\infty \log(1+\lambda_i).
\end{align}
\end{subequations} 

As noted before $m^*$ coincides with the deterministic inverse problem. While this closed-form expression for EIG is only exact for linear inverse problems, it can be utilized as an approximation for nonlinear problems. 

\subsection{Notes on the prior and the function space limit}

In the preceding section we focused on quadratic-norm induced Gaussian measures, so we thereby restrict our attention to the $\mathcal{R}_{L^2}$ and $\mathcal{R}_{H^1}$ regularization utilized in the $\infty-$IDIC for simplicity. We note that there has been work on formulating 
infinite-dimensionally consistent TV priors, and the discussion presented here can be naturally extended to that setting \cite{yao2016tv,lv2020nonlocal}. 

The $L^2(\Omega)$ Tikhonov regularization fails to induce an infinite-dimensionally consistent prior as it is not trace class, since it is induced by the operator $I_\mathcal{M}$
\begin{equation}
	\text{tr}(I_\mathcal{M}) = \sum_{i=1}^\infty \langle \psi_i,\psi_i\rangle_{\mathcal{M}} = \sum_{i=1}^\infty 1 \nless \infty.
\end{equation}

The $H^1(\Omega)$ regularization is associated with the operator $\mathcal{C}^{-1} = -\gamma \Delta + \delta \mathcal{I}_{\mathcal{M}}$ (with homogeneous Neumann boundary conditions), and is only trace class when $\Omega \subset \mathbb{R}^{d}$ for $d=1$. Therefore, the Bayesian extension of the $\mathcal{R}_{H^1}$ problem that we introduced, does not exist in the function space limit. It does however, still lead to a meaningful Bayesian extension of an $\infty-$IDIC problem \emph{with a fixed mesh representation} as in the numerical results that we discuss in section \ref{section:numerics_eig}.

In general one can consider norms based on the $H^k(\Omega)$ semi-norms that do induce meaningful infinite-dimensional Bayesian inverse problems. These are related to the Mat\'{e}rn covariances that are often used in infinite-dimensional Bayesian inverse problems. Particularly, we used them when generating the speckle pattern (see section \ref{appendix:a}) in addition to generating true modulus fields (i.e., the Gaussian random field in Figure \ref{fig:inverse_results_linear}). The differential operator 
\begin{equation}\label{eq:matern_covariance}
	\mathcal{A} = \delta I_\mathcal{M} - \gamma \nabla \cdot (\Theta \nabla),
\end{equation}
can induce a trace class covariance $\mathcal{C} = \mathcal{A}^{-\alpha}$ (with appropriate boundary conditions), when $\alpha > d/2$ \cite{LindgrenRueLindstrom2011}. In \eqref{eq:matern_covariance} $\delta,\gamma\in \mathbb{R}$ control the marginal variance and correlation length in the random fields drawn from the associated Gaussian measure,  and $\Theta \in \mathbb{R}^{d\times d}$ induces anisotropy in the fields \cite{VillaOLearyRoseberry2024}. Consequently, we can take $\delta =0$, and $\Theta = I$ to recover the $H^k(\Omega)$ semi-norms with $\alpha=k$ and see that we require $k>1$ for $\Omega \subset \mathbb{R}^{2}$, but that $k =2$ works for all $d >1$. This enforces even more smoothness in the solutions that are favored by the inverse problems (either Bayesian or the analogous deterministic problem). We do not consider this case in numerical experiments, but note that it leads to a meaningful infinite-dimensional Bayesian analogue of the corresponding deterministic $\infty-$IDIC problem.

\end{document}